\documentclass[12pt]{article}
\usepackage{amssymb}
\usepackage{latexsym}
\usepackage{amsthm}

\newtheorem{definition}{Definition}[section]
\newtheorem{theorem}[definition]{Theorem}
\newtheorem{lemma}[definition]{Lemma}

\newtheorem{remark}[definition]{Remark}

\newtheorem{note}[definition]{Note}

\newtheorem{proposition}[definition]{Proposition}

\newtheorem{notation}[definition]{Notation}

\newcommand{\beast}{\begin{eqnarray*}}
\newcommand{\eeast}{\end{eqnarray*}}

\def\K{\mathbb K}

\def\K{\mathbb K}
\def\fld{\mathbb K}
\def\Mdf{{\hbox{Mat}}_{d+1}(\K)}

\begin{document}

\title{\bf Thin Hessenberg Pairs}
\author{Ali Godjali}


\date{Nov 19, 2009}


\maketitle
\begin{abstract}
\noindent 
A square matrix is called {\it Hessenberg} whenever each entry below the subdiagonal is zero and each entry on the subdiagonal is nonzero. Let $V$ denote a nonzero finite-dimensional vector space over a field $\fld$. 
We consider an ordered pair of linear transformations $A: V
\rightarrow V$ and $A^*: V \rightarrow V$ which satisfy both (i), (ii) below.
\begin{enumerate}
\item There exists a basis for $V$ with respect to which the matrix representing $A$ is Hessenberg and the matrix representing $A^*$ is diagonal. 

\item There exists a basis for $V$ with respect to which the matrix representing $A$ is diagonal and the matrix representing $A^*$ is Hessenberg.  
\end{enumerate}

\noindent We call such a pair a {\it thin Hessenberg pair} (or {\it TH pair}). This is a special case of a {\it Hessenberg pair} which was introduced by the author in an earlier paper. We investigate several bases for $V$ with respect to which the matrices representing $A$ and $A^*$ are attractive.
We display these matrices along with the transition matrices relating the bases. We introduce an ``oriented" version of $A,A^*$
called a TH system. We classify the TH systems up to isomorphism.

\bigskip

\noindent {\bf Keywords}:
Leonard pair, Hessenberg pair, tridiagonal pair, $q$-inverting pair. 
\hfil\break
\noindent {\bf 2010 Mathematics Subject Classification}: 15A04.

\end{abstract}

\section{Introduction}
This paper is about a linear algebraic object called a {\it Hessenberg pair}, which is a generalization of a {\it tridiagonal pair} \cite{TD00, TDRACAH}. We introduced Hessenberg pairs in \cite{hess}. In the present paper, we continue our investigation of Hessenberg pairs by focusing on a special case said to be {\it thin}. To define this case, we will use the following term. A square matrix is called {\it Hessenberg} whenever each entry below the subdiagonal is zero and each entry on the subdiagonal is nonzero. Throughout the paper, $\fld$ will denote a field.

\begin{definition}
\label{def:thinhess} 
\rm 
Let $V$ denote a nonzero finite-dimensional vector space over $\fld$. By a \emph{thin Hessenberg pair} (or {\it TH pair}) on $V$, we mean an ordered pair of linear transformations $A:V\rightarrow V$ and $A^*:V\rightarrow V$ which
satisfy both (i), (ii) below. 

\begin{enumerate}

\item There exists a basis for $V$ with respect to which the matrix representing $A$ is Hessenberg and the matrix representing $A^*$ is diagonal. 

\item There exists a basis for $V$ with respect to which the matrix representing $A$ is diagonal and the matrix representing $A^*$ is Hessenberg.  
\end{enumerate}
\end{definition}

\begin{note}  \samepage
\rm
It is a common notational convention to use $A^*$ to represent the conjugate-transpose of $A$. We are not using this convention. In a TH pair $A, A^*$ the linear transformations $A$ and
$A^*$ are arbitrary subject to (i), (ii) above.
\end{note}

We now briefly summarize the paper. Let $A, A^*$ denote a TH pair on $V$. We investigate six bases for $V$ with respect to which the matrices representing $A$ and $A^*$ are attractive.
We display these matrices along with the transition matrices relating the bases. We introduce an ``oriented" version of $A,A^*$
called a {\it TH system}. We classify the TH systems up to isomorphism. We will give a more detailed summary at the end of Section 2, after establishing some notations and reviewing some basic concepts. 

\section{TH systems}
In our study of a TH pair, it is often helpful to consider a closely related object called a TH system. Before defining this notion,  
we make some definitions and observations. For the rest of the paper, fix an integer $d \geq 0$ and a vector space $V$ over $\fld$ with dimension $d+1$.
Let ${\rm End}(V)$ denote the $\K$-algebra of all linear transformations from $V$ to $V$. Let $\Mdf$ denote the $\K$-algebra consisting of all $(d+1)\times(d+1)$ matrices which have entries in $\K$. We index the rows and columns by $0,1,\ldots, d$. 
Suppose that $\lbrace v_i \rbrace_{i=0}^d$ is a basis for $V$. For $X \in \hbox{Mat}_{d+1}(\K)$ and $Y \in {\rm End}(V)$, we say $X$ \emph{represents} $Y$ \emph{with respect to} $\lbrace v_i \rbrace_{i=0}^d$ whenever $Y v_j = \sum_{i=0}^d X_{ij} v_i$ for $0 \leq j \leq d$. 
For $A \in {\rm End}(V)$ and $W \subseteq V$,  we call $W$ an {\it eigenspace} of
$A$ whenever $W\not=0$ and there exists $\theta \in \K$ such that $W=\lbrace v \in V \;\vert \;Av = \theta v\rbrace$.
In this case $\theta$ is called the {\it eigenvalue} of $A$ corresponding to $W$. We say $A$ is {\it diagonalizable} whenever $V$ is spanned by the eigenspaces of $A$.  We say $A$ is {\it multiplicity-free} whenever $A$ is diagonalizable and each eigenspace of $A$ has dimension one. 

\begin{lemma}
\label{lem:multfree}
Let $A, A^*$ denote a TH pair on $V$. Then each of $A, A^*$ is multiplicity-free. 
\end{lemma}

\noindent {\it Proof:} Concerning $A$, recall by Definition \ref{def:thinhess}(ii) that there exists
a basis for $V$ consisting of eigenvectors of $A$. Consequently, the eigenvalues of $A$ are all in $\fld$ and the minimal polynomial of $A$ has no repeated roots. To show that the eigenvalues of $A$ are distinct, we show that the minimal polynomial of $A$ has degree $d+1$. 
By Definition \ref{def:thinhess}(i), there exists a basis for $V$ with respect to which the matrix representing $A$ is Hessenberg. Denote this matrix by $B$. On one hand, $A$ and $B$ have the same minimal polynomial. On the other hand, using the Hessenberg shape of $B$, we find  $I, B, B^2, \ldots, B^d$ are linearly independent, so the minimal  polynomial of $B$ has degree $d+1$. 
We conclude the minimal polynomial of $A$ has degree $d+1$, so the eigenvalues of $A$ are distinct. Consequently $A$ is multiplicity-free. The case of $A^*$ is similar. \hfill $\Box$

\medskip

\noindent We recall a few more concepts from linear algebra. 
Let $A$ denote a multiplicity-free element of ${\rm End}(V)$. Let $\{V_i\}_{i=0}^d$ denote an ordering of the eigenspaces of $A$ and let $\{\theta_i\}_{i=0}^d$ denote the corresponding ordering of the eigenvalues of $A$. For $0 \leq i \leq d$, define $E_i \in {\rm End}(V)$ 
such that $(E_i-I)V_i=0$ and $E_iV_j=0$ for $j \neq i$ $(0 \leq j \leq d)$. Here $I$ denotes the identity of ${\rm End}(V)$.
We call $E_i$ the {\it primitive idempotent} of $A$ corresponding to $V_i$ (or $\theta_i$).
Observe that
(i)  $I = \sum_{i=0}^d E_i$;
(ii) $E_iE_j=\delta_{i,j}E_i$ $(0 \leq i,j \leq d)$;
(iii) $V_i=E_iV$ $(0 \leq i \leq d)$;
(iv) $A=\sum_{i=0}^d \theta_i E_i$.
Moreover
\begin{equation}        
\label{eq:defEi}
E_i=\prod_{\stackrel{0 \leq j \leq d}{j \neq i}}
          \frac{A-\theta_jI}{\theta_i-\theta_j} \qquad \qquad \qquad (0 \leq i \leq d).
\end{equation}
Note that each of $\{A^i\}_{i=0}^d$, $\{E_i\}_{i=0}^d$ is a basis for the $\K$-subalgebra of $\mbox{\rm End}(V)$ generated by $A$. Moreover $\prod_{i=0}^d(A-\theta_iI)=0$.

\medskip
\noindent We now define a TH system. 

\begin{definition} 
\label{def:HS}
\rm
By a {\it thin Hessenberg system} (or {\it TH system}) on $V$ we mean a sequence
\[
\Phi=(A;\lbrace E_i\rbrace_{i=0}^d;A^*;\lbrace E^*_i\rbrace_{i=0}^d)
\]
which satisfies (i)--(v) below. 
\begin{enumerate}
\item 
Each of $A,A^*$ is a multiplicity-free element of ${\rm End}(V)$.
\item 
$\{E_i\}_{i=0}^d$ is an ordering of the primitive idempotents of $A$.
\item 
$\{E^*_i\}_{i=0}^d$ is an ordering of the primitive idempotents of $A^*$.
\item 
${\displaystyle{
E_iA^*E_j = \cases{0, &if $\; i-j > 1$\cr
\not=0, &if $\; i-j  = 1$\cr}
\qquad \qquad 
(0 \leq i,j\leq d)}}$.
\item 
${\displaystyle{
 E^*_iAE^*_j = \cases{0, &if $\; i-j > 1$\cr
\not=0, &if $\; i-j  = 1$\cr}
\qquad \qquad 
(0 \leq i,j\leq d).}}$
\end{enumerate}
We call $V$ the {\it underlying} vector space and say $\Phi$ is {\it over} $\K$.  
\end{definition}

\noindent We comment on how TH pairs and TH systems are related.
Let $(A;\lbrace E_i\rbrace_{i=0}^d;A^*;\lbrace E^*_i\rbrace_{i=0}^d)$ denote a TH system on $V$. 
For $0 \leq i \leq d$, let $v_i$ (resp. $v^*_i$) denote a nonzero vector in $E_iV$ (resp. $E^*_iV$). Then the sequence $\lbrace v_i \rbrace_{i=0}^d$ (resp. $\lbrace v^*_i \rbrace_{i=0}^d$) is a basis for $V$ which satisfies Definition \ref{def:thinhess}(ii) (resp. Definition \ref{def:thinhess}(i)). Therefore the pair $A, A^*$ is a TH pair on $V$. Conversely, let $A, A^*$ denote a TH pair on $V$. Then each of $A, A^*$ is multiplicity-free by Lemma \ref{lem:multfree}. Let $\lbrace v_i \rbrace_{i=0}^d$ (resp. $\lbrace v^*_i \rbrace_{i=0}^d$) denote a basis for $V$ which satisfies Definition \ref{def:thinhess}(ii) (resp. Definition \ref{def:thinhess}(i)). For $0 \leq i \leq d$, the vector $v_i$ (resp. $v^*_i$) is an eigenvector for $A$ (resp. $A^*$); let $E_i$ (resp. $E^*_i$) denote the corresponding primitive idempotent. Then $(A;\lbrace E_i\rbrace_{i=0}^d;A^*;\lbrace E^*_i\rbrace_{i=0}^d)$ is a TH system on $V$. 

\begin{definition}
\label{def:thpairassthsytem}
\rm
Let $\Phi = (A;\lbrace E_i\rbrace_{i=0}^d;A^*;\lbrace E^*_i\rbrace_{i=0}^d)$ denote a TH system on $V$.
Observe that $A, A^*$ is a TH pair on $V$. We say this pair is {\it associated} with $\Phi$. 
\end{definition}

\begin{remark}
\rm
With reference to Definition \ref{def:thpairassthsytem}, conceivably a given TH pair is associated with many TH systems. 
\end{remark}

\medskip
\noindent We now define the notion of {\it isomorphism} for TH systems. 

\begin{definition}
\rm
Let $\Phi=(A;\lbrace E_i\rbrace_{i=0}^d;A^*;\lbrace E^*_i\rbrace_{i=0}^d)$ denote a TH system on $V$. Let $W$ denote a vector space over $\fld$ with dimension $d+1$, and let $\Psi:=(B;\lbrace F_i \rbrace_{i=0}^d;B^*;\lbrace F^*_i\rbrace_{i=0}^d)$ denote a TH system on $W$. By an {\it isomorphism of TH systems} from $\Phi$ to $\Psi$ we mean a vector space isomorphism $\gamma: V \rightarrow W$ such that 
\begin{eqnarray*}
\label{eq:isosys}
\gamma A
=B\gamma,
\qquad 
\gamma A^*=
B^*\gamma, \qquad \gamma E_i = F_i \gamma \ \ \   (0 \leq i \leq d), \qquad \gamma E^*_i = F^*_i \gamma \ \ \ (0 \leq i \leq d).
\end{eqnarray*}
\noindent We say $\Phi$ and $\Psi$ are {\it isomorphic} whenever there exists an isomorphism of TH systems from $\Phi$ to $\Psi$. 
\end{definition}

\noindent We now define the dual of a TH system. 

\begin{definition}
\label{def:THdual}
\rm
Let $\Phi=(A;\lbrace E_i\rbrace_{i=0}^d;A^*;\lbrace E^*_i\rbrace_{i=0}^d)$ denote a TH system on $V$. Observe that 
$\Phi^*:=(A^*;\lbrace E^*_i\rbrace_{i=0}^d;A;\lbrace E_i\rbrace_{i=0}^d)$ is a TH system on $V$. We call $\Phi^*$ the {\it dual} of $\Phi$. 
\end{definition}

We recall some more terms. Let $\lbrace v_i \rbrace_{i=0}^d$ denote a basis for $V$. By the {\it inversion} of $\lbrace v_i \rbrace_{i=0}^d$ we mean the basis $\lbrace v_{d-i} \rbrace_{i=0}^d$. A square matrix is called {\it lower bidiagonal} whenever each nonzero entry lies on either the diagonal or the subdiagonal, and each entry on the subdiagonal is nonzero. A square matrix is called {\it upper bidiagonal} whenever its transpose is lower bidiagonal. 

We now give a detailed summary of the paper.
Let $\Phi=(A;\lbrace E_i\rbrace_{i=0}^d;A^*;\lbrace E^*_i\rbrace_{i=0}^d)$ denote a TH system on $V$. We investigate six bases for $V$ that we find attractive. The first four are called the {\it $\Phi$-split basis}, the {\it $\Phi^*$-split basis}, the {\it inverted $\Phi$-split basis}, and the {\it inverted $\Phi^*$-split basis}. With respect to each of these bases, the matrix representing one of $A, A^*$ is lower bidiagonal and the matrix representing the other is upper bidiagonal. The other two bases in our investigation are called the {\it $\Phi$-standard basis} and the {\it $\Phi^*$-standard basis}. A $\Phi$-standard basis (resp. $\Phi^*$-standard basis) satisfies Definition \ref{def:thinhess}(i) (resp. Definition \ref{def:thinhess}(ii)), subject to a certain normalization. For each of the six bases, we display the matrices representing $A$ and $A^*$.  
We display some transition matrices relating these six bases. We associate with $\Phi$ a sequence of scalars called its {\it parameter array}. We show that  $\Phi$ is determined up to isomorphism by its parameter array. Using this fact, we classify the TH systems up to isomorphism.  

\section{The eigenvalue sequences}
Let $\Phi$ denote a TH system. In this section we associate with $\Phi$ two sequences of scalars called the {\it eigenvalue sequence} and the {\it dual eigenvalue sequence}. We describe some properties of these sequences that we will use later in the paper.

\begin{definition}
\label{def:evseq}
\rm
Let $\Phi=(A;\lbrace E_i\rbrace_{i=0}^d;A^*;\lbrace E^*_i\rbrace_{i=0}^d)$ denote a TH system on $V$.
For $0 \leq i \leq d$, let $\theta_i $ (resp. $\theta^*_i$) denote the eigenvalue of $A$ (resp. $A^*$) corresponding to $E_i$ (resp. $E^*_i$). We refer to $\lbrace \theta_i \rbrace_{i=0}^d$ as the {\it eigenvalue sequence} of $\Phi$. We refer to  $\lbrace \theta^*_i \rbrace_{i=0}^d$ as the {\it dual eigenvalue sequence} of $\Phi$. We observe that $\lbrace \theta_i \rbrace_{i=0}^d$ are mutually distinct and contained in $\K$. Similarly $\lbrace \theta^*_i \rbrace_{i=0}^d$ are mutually distinct and contained in $\K$. 
\end{definition}

\begin{definition}
\label{def:evseqhp}
\rm 
Let $A, A^*$ denote a TH pair. By an {\it eigenvalue sequence} of $A, A^*$, we mean the eigenvalue sequence of an associated TH system. 
By a {\it dual eigenvalue sequence} of $A, A^*$, we mean an eigenvalue sequence of the TH pair $A^*, A$. We emphasize that a given TH pair could have many eigenvalue and dual eigenvalue sequences.
\end{definition}
 
\begin{lemma}
\label{lem:eig}
Let $A, A^*$ denote a TH pair. Then the following {\rm (i), (ii)} hold.

\begin{enumerate}
\item[\rm (i)] 
Let $\lbrace \theta_i \rbrace_{i=0}^d$ and $\lbrace \theta'_i \rbrace_{i=0}^d$ denote eigenvalue sequences of $A, A^*$ such that $\theta_0 = \theta'_0$. Then $\theta_i  = \theta'_i$ for $0 \leq i \leq d$.
\item[\rm (ii)]
Let $\lbrace \theta^*_i \rbrace_{i=0}^d$ and $\lbrace \theta^{*\prime}_i \rbrace_{i=0}^d$ denote dual eigenvalue sequences of $A, A^*$ such that $\theta^*_0 = \theta^{*\prime}_0$. Then $\theta^*_i = \theta^{*\prime}_i$ for $0 \leq i \leq d$.
\end{enumerate}
\end{lemma}

\noindent {\it Proof:} (i)
For $0 \leq i \leq d$, let $V_i$ (resp. $V_i'$) denote the eigenspace of $A$ corresponding to $\theta_i$ (resp. $\theta'_i$). 
It suffices to show that $V_{i} = V_{i}'$ for $0 \leq i \leq d$. This follows once we show that $W_i = W_i'$ for $0 \leq i \leq d$, where $W_i = \sum_{h=0}^i V_h$ and $W_i' = \sum_{h=0}^i V_h'$. We prove this by induction on $i$. First assume that $i=0$. Then $W_0 = V_0 = V_0' = W_0'$ since $\theta_0 = \theta'_0$. Next assume that $1 \leq i \leq d$. By induction, we have $W_{i-1} = W_{i-1}'$. By Definition \ref{def:thinhess}(ii) and since $\lbrace \theta_h \rbrace_{h=0}^d$ is an eigenvalue sequence of $A, A^*$,  we find $A^* W_{i-1} \subseteq W_{i}$ and $A^* W_{i-1} \nsubseteq W_{i-1}$. Therefore $W_{i-1} + A^* W_{i-1} =  W_{i}$. Similarly  $W_{i-1}' + A^* W_{i-1}' =  W_{i}'$.  Comparing these equations using $W_{i-1} = W_{i-1}'$, we find $W_i = W_i'$. The result follows. 

\smallskip
\noindent (ii) Apply (i) to $\Phi^*$. \hfill $\Box$

\begin{lemma}
\label{lem:trihess}
Let $A, A^*$ denote a TH pair on $V$. Fix a basis for $V$ and let $B$ (resp. $B^*$) denote the matrix in $\hbox{Mat}_{d+1}(\K)$ which represents $A$ (resp. $A^*$) with respect to that basis. Assume that $B$ is Hessenberg and $B^*$ is upper triangular. Then the sequence of diagonal entries $\lbrace B^*_{ii} \rbrace_{i=0}^d$ of $B^*$ is a dual eigenvalue sequence of $A, A^*$.
\end{lemma}
\noindent {\it Proof:}
We assume that $A, A^*$ is a TH pair so $A^*$ is multiplicity-free. We assume that $B^*$ is upper triangular so the sequence $\lbrace B^*_{ii} \rbrace_{i=0}^d$ is an ordering of the eigenvalues of $A^*$. We show that this sequence is a dual eigenvalue sequence of $A, A^*$. For $0 \leq i \leq d$, let $E^*_i$ denote the primitive idempotent of $A^*$ corresponding to the eigenvalue $B^*_{ii}$. It suffices to show that Definition \ref{def:HS}(v) holds. We make a few observations. Let $\lbrace v_i \rbrace_{i=0}^d$ denote the basis for $V$ in the statement of the lemma. For $0 \leq i \leq d$, let $W_i$ denote the subspace of $V$ spanned by $\lbrace v_h \rbrace_{h=0}^i$. The matrix $B^*$ is upper triangular so
$A^*W_i \subseteq W_i$. The restriction of $A^*$ to $W_i$ has eigenvalues $\lbrace B^*_{hh} \rbrace_{h=0}^i$, so
$W_i=E^*_0V+\cdots + E^*_iV$. 
Since $B$ is Hessenberg, we have $A E^*_jV \subseteq W_{j+1}$ and $A E^*_jV \nsubseteq W_j$ for $0 \leq j \leq d-1$.  
We can now easily show that Definition \ref{def:HS}(v) holds.
Fix integers $i,j$ $(0 \leq i,j \leq d)$ such that $i-j \geq 1$. First assume that $i-j>1$.
From our observations above, $E^*_i A E^*_j V \subseteq E^*_i W_{j+1} = 0$.
Therefore $E^*_i A E^*_j V =0$ so $E^*_i A E^*_j =0$.
Next assume that $i-j=1$. We show $E^*_i A E^*_j \not=0$. By way of contradiction, assume that $E^*_i A E^*_j =0$.
By this and our earlier observations, we have $E^*_h A E^*_j=0$ for $i \leq h \leq d$.
Therefore $AE^*_j = \sum_{h=0}^d E^*_h A E^*_j = \sum_{h=0}^j E^*_h A E^*_j$, so
$AE^*_jV \subseteq \sum_{h=0}^j E^*_h A E^*_jV \subseteq W_j$. This contradicts our above remarks so
$E^*_i A E^*_j \not=0$. The result follows. \hfill $\Box$

\section{The $\Phi$-split basis}
Let $\Phi$ denote a TH system on $V$. In this section we investigate a certain basis for $V$ called the $\Phi$-split basis. We will refer to the following notation. 

\begin{notation}
\label{not:aastar}
\rm
Let $\Phi=(A;\lbrace E_i\rbrace_{i=0}^d;A^*;\lbrace E^*_i\rbrace_{i=0}^d)$ denote a TH system on $V$. Let $\lbrace \theta_i \rbrace_{i=0}^d$ (resp. $\lbrace \theta^*_i \rbrace_{i=0}^d$) denote the eigenvalue (resp. dual eigenvalue) sequence of $\Phi$. 
\end{notation}

With reference to Notation \ref{not:aastar}, in our study of $\Phi$ we will use the following term. 
By a {\it decomposition} of $V$ we mean a sequence $\lbrace U_i \rbrace_{i=0}^d$ of one-dimensional subspaces of $V$ such that
\begin{eqnarray*}
\qquad \qquad V=U_0+U_1+\cdots + U_d \qquad
\qquad (\hbox{direct sum}).
\end{eqnarray*}
For notational convenience, set $U_{-1}=0$ and $U_{d+1}=0$. We now describe a certain decomposition of $V$ associated with $\Phi$. For $0 \leq i \leq d$, define 

\begin{equation}
\label{eq:defui}
U_i = (E^*_0V + E^*_1V + \cdots + E^*_iV) \cap (E_0V + E_{1}V + \cdots + E_{d-i}V).
\end{equation}

\noindent By \cite[Lemma 2.4, Theorem 4.1]{split}, the sequence $\lbrace U_i \rbrace_{i=0}^d$ is a decomposition of $V$. Moreover for $0 \leq i \leq d$, both
\begin{eqnarray}
(A-\theta_{d-i} I)U_i &=& U_{i+1},
\label{eq:raise}
\\
(A^*-\theta^*_i I)U_i &=& U_{i-1}.
\label{eq:lower}
\end{eqnarray}

\noindent Setting $i=d$ in (\ref{eq:defui}) we find $U_d=E_0V$. Combining this with
(\ref{eq:lower}) we find  
\begin{equation}
U_i = (A^* -\theta^*_{i+1} I)\cdots (A^* - \theta^*_{d-1} I)(A^* -\theta^*_d I)E_0V
\qquad \qquad (0 \leq i\leq d).
\label{eq:uialt}
\end{equation}
Let $\eta_0$ denote a nonzero vector in $E_0V$. From (\ref{eq:uialt}) we find that for 
$0 \leq i \leq d$, the vector $(A^*-\theta^*_{i+1}I)\cdots (A^*-\theta^*_{d}I)\eta_0 $ is a basis for
$U_i$. By this and since $\lbrace U_i \rbrace_{i=0}^d$ is a decomposition of $V$, the sequence
\begin{eqnarray*}
(A^*-\theta^*_{i+1} I)\cdots (A^*-\theta^*_{d-1} I)(A^*-\theta^*_{d}I)\eta_0  \qquad \qquad
(0 \leq i \leq d)
\end{eqnarray*}
is a basis for $V$.

\begin{definition}
\label{def:lbubbasis}
\rm
With reference to Notation \ref{not:aastar}, a basis for $V$ is said to be {\it $\Phi$-split} whenever it is of the form 
\begin{equation}
\label{eq:spbasis}
(A^*-\theta^*_{i+1} I)\cdots (A^*-\theta^*_{d-1} I)(A^*-\theta^*_{d}I)\eta_0  \qquad \qquad (0 \leq i \leq d),
\end{equation}
where $0 \neq \eta_0 \in E_0V$. 
\end{definition}

\begin{lemma}
\label{lem:scalar}
With reference to Notation \ref{not:aastar}, let $\lbrace v_i \rbrace_{i=0}^d$ denote a $\Phi$-split basis for $V$, and let $\lbrace w_i \rbrace_{i=0}^d$ denote any vectors in $V$. Then the following are equivalent. 

\begin{enumerate}
\item[\rm (i)] 
$\lbrace w_i \rbrace_{i=0}^d$ is a $\Phi$-split basis for $V$. 
\item[\rm (ii)]
There exists $0 \neq c \in \K$ such that $w_i = c\,v_i$ for $0 \leq i \leq d$.  
\end{enumerate}
\end{lemma}

\noindent {\it Proof:} Routine. \hfill $\Box$

\medskip

\noindent With reference to Notation \ref{not:aastar}, our next goal is to describe the matrices representing $A, A^*$ with respect to a $\Phi$-split basis for $V$. We start with an observation. 
Let $1 \leq i \leq d$. By (\ref{eq:lower}) we have $(A^*-\theta^*_i I)U_i = U_{i-1}$, and by (\ref{eq:raise}) we have $(A-\theta_{d-i+1} I)U_{i-1} = U_i$. Therefore $U_i$ is an eigenspace of $(A-\theta_{d-i+1}I)(A^*-\theta^*_i I)$ and the corresponding eigenvalue is a nonzero element of $\K$. We denote this eigenvalue by $\phi_i$. We call the sequence $\lbrace \phi_i \rbrace_{i=1}^d$ the {\it split sequence} of $\Phi$. For notational convenience, set $\phi_0=0$ and $\phi_{d+1} = 0$.

\begin{proposition}
\label{prop:lbublooklike}
With reference to Notation \ref{not:aastar}, let $B$ (resp. $B^*$) denote the matrix in $\hbox{Mat}_{d+1}(\K)$ which represents $A$ (resp. $A^*$) with respect to a $\Phi$-split basis for $V$. Then 
\begin{equation}
\label{eq:matrepaastar}
B = 
\left(
\begin{array}{c c c c c c}
\theta_d & & & & & {\bf 0} \\
\phi_1 & \theta_{d-1} &  & & & \\
& \phi_2 & \theta_{d-2} &  & & \\
& & \cdot & \cdot &  &  \\
& & & \cdot & \cdot &  \\
{\bf 0}& & & & \phi_d & \theta_0
\end{array}
\right),
\qquad \quad 
B^* = 
\left(
\begin{array}{c c c c c c}
\theta^*_0 & 1 & & & & {\bf 0} \\
& \theta^*_1 & 1 & & & \\
& & \theta^*_2 & \cdot & & \\
& &  & \cdot & \cdot &  \\
& & &  & \cdot & 1 \\
{\bf 0}& & & &  & \theta^*_d
\end{array}
\right),
\end{equation}
\noindent where $\lbrace \phi_i \rbrace_{i=1}^d$ is the split sequence of $\Phi$. In particular, $B$ is lower bidiagonal and $B^*$ is upper bidiagonal. 
\end{proposition}
\noindent {\it Proof:} Follows from Definition \ref{def:lbubbasis} and the discussion prior to this proposition. \hfill $\Box$

\bigskip

\noindent We give an alternate description of the $\Phi$-split basis. To motivate this, we set $i=0$ in (\ref{eq:defui}) and find $U_0=E^*_0V$. Combining this with (\ref{eq:raise}) we find  
\begin{equation}
\label{eq:uialtvar}
U_i = (A -\theta_{d-i+1} I)\cdots (A - \theta_{d-1} I)(A -\theta_d I)E^*_0V
\qquad \qquad (0 \leq i\leq d).
\end{equation}
Let $\eta^*_0$ denote a nonzero vector in $E^*_0V$. From (\ref{eq:uialtvar}) we find that for 
$0 \leq i \leq d$, the vector $(A -\theta_{d-i+1} I) \cdots (A -\theta_d I)\eta^*_0 $ is a basis for
$U_i$. By this and since $\lbrace U_i \rbrace_{i=0}^d$ is a decomposition of $V$, the sequence
\begin{eqnarray*}
(A -\theta_{d-i+1} I)\cdots (A - \theta_{d-1} I)(A -\theta_d I)\eta^*_0  \qquad \qquad
(0 \leq i \leq d)
\end{eqnarray*}
is a basis for $V$. This basis is not a $\Phi$-split basis in general, but we do have the following result. 

\begin{lemma}
\label{lem:lbubbasisvar}
With reference to Notation \ref{not:aastar}, let $\lbrace v_i \rbrace_{i=0}^d$ denote any vectors in $V$. Then the following are equivalent.

\begin{enumerate}
 \item[\rm (i)]
$\lbrace v_i \rbrace_{i=0}^d$ is a $\Phi$-split basis for $V$.
 \item[\rm (ii)]
There exists $0 \neq \eta^*_0 \in E^*_0V$ such that
\begin{equation}
\label{eq:lbubbasisvar}
\displaystyle{
v_i = \frac{(A-\theta_{d-i+1} I)\cdots (A-\theta_{d-1} I)(A-\theta_{d}I)\eta^*_0}{\phi_1 \phi_2 \cdots \phi_i} \qquad \qquad (0 \leq i \leq d),}
\end{equation}
where $\lbrace \phi_i \rbrace_{i=1}^d$ is the split sequence of $\Phi$.
\end{enumerate}
\end{lemma}

\noindent {\it Proof:} Suppose that $\lbrace v_i \rbrace_{i=0}^d$ is a $\Phi$-split basis for $V$. Then by Proposition \ref{prop:lbublooklike}, the matrices representing $A$ and $A^*$ with respect to $\lbrace v_i \rbrace_{i=0}^d$ are as shown in (\ref{eq:matrepaastar}). From the matrix on the right in (\ref{eq:matrepaastar}) we find $v_0 \in E^*_0V$. Taking $\eta^*_0 = v_0$, we find (\ref{eq:lbubbasisvar}) holds for $i=0$. From the matrix on the left in (\ref{eq:matrepaastar}), we find that (\ref{eq:lbubbasisvar}) holds for $1 \leq i \leq d$. The result follows by Lemma \ref{lem:scalar}. \hfill $\Box$

\medskip

\noindent We now give an alternate description of the split sequence of $\Phi$. 

\begin{lemma}
\label{lem:splitseqvar}
With reference to Notation \ref{not:aastar}, let $\lbrace \phi_i \rbrace_{i=1}^d$ denote the split sequence of $\Phi$. Then for $1 \leq i \leq d$, $\phi_i$ is the eigenvalue of $(A^* - \theta^*_{i}I)(A - \theta_{d-i+1}I)$ for the eigenspace $U_{i-1}$, where $\lbrace U_i \rbrace_{i=0}^d$ is from (\ref{eq:defui}).
\end{lemma}

\noindent {\it Proof:} Let $\lbrace v_j \rbrace_{j=0}^d$ denote a $\Phi$-split basis for $V$. Recall that $v_j$ spans $U_j$ for $0 \leq j \leq d$, so it suffices to show that $(A^* - \theta^*_{i}I)(A - \theta_{d-i+1}I) v_{i-1} =  \phi_i v_{i-1}$. This follows by Proposition \ref{prop:lbublooklike} and a routine matrix computation. \hfill $\Box$

\medskip

\noindent We now give three characterizations of the $\Phi$-split basis. 

\begin{proposition}
\label{lem:lbubbasis}
With reference to Notation \ref{not:aastar}, let $\lbrace v_i \rbrace_{i=0}^d$ denote a sequence of vectors in $V$, not all zero. Then this sequence is a $\Phi$-split basis for $V$ if and only if the following {\rm (i), (ii)} hold.
\begin{enumerate}
\item[\rm (i)] 
$v_d \in E_0V$.
\item[\rm (ii)]
$A^*v_i = \theta^*_i v_i + v_{i-1}$ for $1\leq i \leq d.$ 
\end{enumerate}
\end{proposition}

\noindent {\it Proof:} Routine using Definition \ref{def:lbubbasis}. \hfill $\Box$

\begin{proposition}
\label{lem:lbubbasis2}
With reference to Notation \ref{not:aastar}, let $\lbrace v_i \rbrace_{i=0}^d$ denote a sequence of vectors in $V$, not all zero. Then this sequence is a $\Phi$-split basis for $V$ if and only if the following {\rm (i), (ii)} hold.
\begin{enumerate}
\item[\rm (i)] 
$v_0 \in E^*_0V$.
\item[\rm (ii)]
$A v_i = \theta_{d-i} v_i + \phi_{i+1} v_{i+1}$ for $0 \leq i \leq d-1$, where $\lbrace \phi_i \rbrace_{i=1}^d$ is the split sequence of $\Phi$. 
\end{enumerate}
\end{proposition}

\noindent {\it Proof:} Routine using Lemma \ref{lem:lbubbasisvar}. \hfill $\Box$

\begin{proposition}
\label{thm:lbubbasisvsrep}
With reference to Notation \ref{not:aastar}, let $\lbrace v_i \rbrace_{i=0}^d$ denote a basis for $V$, and let $C$ (resp. $C^*$) denote the matrix in $\hbox{Mat}_{d+1}(\K)$ which represents $A$ (resp. $A^*$) with respect to this basis. Then $\lbrace v_i \rbrace_{i=0}^d$ is a $\Phi$-split basis for $V$ if and only if the following {\rm (i)--(iii)} hold.

\begin{enumerate}
\item[\rm (i)]
$C$ is lower bidiagonal and $C^*$ is upper bidiagonal. 
\item[\rm (ii)]
$C^*_{i-1,i} = 1$ for $1 \leq i \leq d$.
\item[\rm (iii)]
$C_{dd} = \theta_0$ and $C^*_{00} = \theta^*_0$.
\end{enumerate}
\end{proposition}

\noindent {\it Proof:} Suppose that $\lbrace v_i \rbrace_{i=0}^d$ is a $\Phi$-split basis for $V$. By Proposition \ref{prop:lbublooklike}, conditions (i)--(iii) above hold. We have proved the proposition in one direction. For the other direction, suppose that conditions (i)--(iii) above hold. To show that $\lbrace v_i \rbrace_{i=0}^d$ is a $\Phi$-split basis for $V$, we invoke Proposition \ref{lem:lbubbasis}. We show that Proposition \ref{lem:lbubbasis}(i), (ii) hold. 
Since $C$ is lower bidiagonal with $C_{dd} = \theta_0$, $v_d$ is an eigenvector of $A$ corresponding to eigenvalue $\theta_0$. Therefore Proposition \ref{lem:lbubbasis}(i) holds. 
Since $C$ is lower bidiagonal, $C$ is Hessenberg. Since $C^*$ is upper bidiagonal, $C^*$ is upper triangular. Now apply Lemma \ref{lem:trihess} to conclude that $\lbrace C^*_{ii} \rbrace_{i=0}^d$ is a dual eigenvalue sequence of $A, A^*$. Therefore since $C^*_{00} = \theta^*_0$, by Lemma \ref{lem:eig} we have $C^*_{ii} = \theta^*_i$ for $0 \leq i \leq d$. Since $C^*$ is upper bidiagonal with $C^*_{i-1,i} = 1$ for $1 \leq i \leq d$, we have $A^*v_i = \theta^*_{i}v_i + v_{i-1}$ for $1\leq i \leq d$. Therefore Proposition \ref{lem:lbubbasis}(ii) holds. The result follows. \hfill $\Box$

\section{Variations on the $\Phi$-split basis}
Let $\Phi$ denote a TH system on $V$. In the previous section we discussed the $\Phi$-split basis for $V$. In this section we discuss three variations on this basis called the $\Phi^*$-split basis, the inverted $\Phi$-split basis, and the inverted $\Phi^*$-split basis. The following lemma will be useful. 

\begin{lemma}
\label{lem:splitseqdual}
With reference to Notation \ref{not:aastar}, let $\lbrace \phi_i \rbrace_{i=1}^d$ denote the split sequence of $\Phi$. Then the split sequence of $\Phi^*$ is $\lbrace \phi_{d-i+1} \rbrace_{i=1}^d$.
\end{lemma}

\noindent {\it Proof:} Routine using the definition of the split sequence of $\Phi$ and Lemma \ref{lem:splitseqvar}. \hfill $\Box$

\medskip

\noindent We now discuss the  $\Phi^*$-split basis.

\begin{lemma}
\label{lem:duallbubbasis}
With reference to Notation \ref{not:aastar}, let $\lbrace v_i \rbrace_{i=0}^d$ denote any vectors in $V$. Then the following are equivalent.

\begin{enumerate}
 \item[\rm (i)]
$\lbrace v_i \rbrace_{i=0}^d$ is a $\Phi^*$-split basis for $V$.
 \item[\rm (ii)]
There exists $0 \neq \eta^*_0 \in E^*_0V$ such that
\begin{equation}
\label{eq:duallbubbasis}
v_i = (A-\theta_{i+1} I)\cdots (A-\theta_{d-1} I)(A-\theta_{d}I)\eta^*_0  \qquad \qquad (0 \leq i \leq d).
\end{equation}
 \item[\rm (iii)]
There exists $0 \neq \eta_0 \in E_0V$ such that
\begin{equation}
\label{eq:duallbubbasisvar}
\displaystyle{
v_i = \frac{(A^*-\theta^*_{d-i+1} I)\cdots (A^*-\theta^*_{d-1} I)(A^*-\theta^*_{d}I)\eta_0}{\phi_{d-i+1} \cdots \phi_{d-1} \phi_d}  \qquad (0 \leq i \leq d),}
\end{equation}
where $\lbrace \phi_i \rbrace_{i=1}^d$ is the split sequence of $\Phi$.
\end{enumerate}
\end{lemma}

\noindent {\it Proof:} To prove (i) $\leftrightarrow$ (ii) apply Definition \ref{def:lbubbasis} to $\Phi^*$.
To prove (i) $\leftrightarrow$ (iii) apply Lemma \ref{lem:lbubbasisvar} to $\Phi^*$ and use Lemma \ref{lem:splitseqdual}. \hfill $\Box$

\medskip

\begin{proposition}
\label{prop:duallbublooklike}
With reference to Notation \ref{not:aastar}, let $B$ (resp. $B^*$) denote the matrix in $\hbox{Mat}_{d+1}(\K)$ which represents $A$ (resp. $A^*$) with respect to a $\Phi^*$-split basis for $V$. Then 
\begin{eqnarray*}
B = 
\left(
\begin{array}{c c c c c c}
\theta_0 & 1 & & & & {\bf 0} \\
& \theta_1 & 1 & & & \\
& & \theta_2 & \cdot & & \\
& & & \cdot & \cdot &  \\
& & & & \cdot & 1 \\
{\bf 0}& & & & & \theta_d
\end{array}
\right), \qquad  \quad 
B^* = 
\left(
\begin{array}{c c c c c c}
\theta^*_d & & & & & {\bf 0} \\
\phi_d & \theta^*_{d-1} & & & & \\
& \phi_{d-1} & \theta^*_{d-2} & & & \\
& & \cdot & \cdot & &  \\
& & & \cdot & \cdot &  \\
{\bf 0}& & & & \phi_1 & \theta^*_0
\end{array}
\right),
\end{eqnarray*}
\noindent where $\lbrace \phi_i \rbrace_{i=1}^d$ is the split sequence of $\Phi$. In particular, $B$ is upper bidiagonal and $B^*$ is lower bidiagonal. 
\end{proposition}

\noindent {\it Proof:} Apply Proposition \ref{prop:lbublooklike} to $\Phi^*$ and use Lemma \ref{lem:splitseqdual}. \hfill $\Box$

\begin{lemma}
\label{cor:trans}
With reference to Notation \ref{not:aastar}, let $\lbrace v_i \rbrace_{i=0}^d$ denote a $\Phi$-split basis for $V$, and let $\lbrace w_i \rbrace_{i=0}^d$ denote any vectors in $V$. Then the following are equivalent. 
\begin{enumerate}
\item[\rm (i)] 
$\lbrace w_i \rbrace_{i=0}^d$ is a $\Phi^*$-split basis for $V$. 
\item[\rm (ii)]
There exists $0 \neq c \in \K$ such that $\displaystyle{w_i = \frac{c \, v_{d-i}}{\phi_{d-i+1} \cdots \phi_{d-1} \phi_d}}$ for $0 \leq i \leq d$,
where $\lbrace \phi_i \rbrace_{i=1}^d$ is the split sequence of $\Phi$.  
\end{enumerate}
\end{lemma}

\noindent {\it Proof:} Compare (\ref{eq:spbasis}) and (\ref{eq:duallbubbasisvar}). \hfill $\Box$

\medskip
\noindent We now turn our attention to the inverted $\Phi$-split basis. 

\begin{lemma}
\label{lem:invertedlbubbasis}
With reference to Notation \ref{not:aastar}, let $\lbrace v_i \rbrace_{i=0}^d$ denote any vectors in $V$. Then the following are equivalent.

\begin{enumerate}
 \item[\rm (i)]  
$\lbrace v_i \rbrace_{i=0}^d$ is an inverted $\Phi$-split basis for $V$.
 \item[\rm (ii)] 
There exists $0 \neq \eta_0 \in E_0V$ such that
\begin{eqnarray*}
v_i = (A^*-\theta^*_{d-i+1} I)\cdots (A^*-\theta^*_{d-1} I)(A^*-\theta^*_{d}I)\eta_0 \qquad \qquad (0 \leq i \leq d).
\end{eqnarray*}
 \item[\rm (iii)] 
There exists $0 \neq \eta^*_0 \in E^*_0V$ such that 
\begin{eqnarray*}
\displaystyle{v_i = \frac{(A-\theta_{i+1} I)\cdots (A-\theta_{d-1} I)(A-\theta_{d}I)\eta^*_0}{\phi_1 \phi_2 \cdots \phi_{d-i}}  \qquad \qquad (0 \leq i \leq d),}
\end{eqnarray*}
where $\lbrace \phi_i \rbrace_{i=1}^d$ is the split sequence of $\Phi$.
\end{enumerate}
\end{lemma}

\noindent {\it Proof:} Routine using Definition \ref{def:lbubbasis}, Lemma \ref{lem:lbubbasisvar}, and the meaning of inversion. \hfill $\Box$

\begin{proposition}
\label{prop:invertedlbublooklike}
With reference to Notation \ref{not:aastar}, let $B$ (resp. $B^*$) denote the matrix in $\hbox{Mat}_{d+1}(\K)$ which represents $A$ (resp. $A^*$) with respect to an inverted $\Phi$-split basis for $V$. Then 
\begin{eqnarray*}
B = 
\left(
\begin{array}{c c c c c c}
\theta_0 & \phi_d & & & & {\bf 0} \\
& \theta_1 & \phi_{d-1} & & & \\
& & \theta_2 & \cdot & & \\
& & & \cdot & \cdot &  \\
& & & & \cdot & \phi_1 \\
{\bf 0}& & & & & \theta_d
\end{array}
\right), \qquad  \quad 
B^* = 
\left(
\begin{array}{c c c c c c}
\theta^*_d & & & & & {\bf 0} \\
1 & \theta^*_{d-1} & & & & \\
& 1 & \theta^*_{d-2} & & & \\
& & \cdot & \cdot & &  \\
& & & \cdot & \cdot &  \\
{\bf 0}& & & & 1 & \theta^*_0
\end{array}
\right),
\end{eqnarray*}
\noindent where $\lbrace \phi_i \rbrace_{i=1}^d$ is the split sequence of $\Phi$. In particular, $B$ is upper bidiagonal and $B^*$ is lower bidiagonal. 
\end{proposition}

\noindent {\it Proof:} Routine using Proposition \ref{prop:lbublooklike} and the meaning of inversion. \hfill $\Box$

\begin{lemma}
With reference to Notation \ref{not:aastar}, let $\lbrace v_i \rbrace_{i=0}^d$ denote a $\Phi$-split basis for $V$, and let $\lbrace w_i \rbrace_{i=0}^d$ denote any vectors in $V$. Then the following are equivalent. 
\begin{enumerate}
\item[\rm (i)] 
$\lbrace w_i \rbrace_{i=0}^d$ is an inverted $\Phi$-split basis for $V$. 
\item[\rm (ii)] 
There exists $0 \neq c \in \K$ such that $w_i = c \,v_{d-i}$ for $0 \leq i \leq d$.
\end{enumerate}
\end{lemma}

\noindent {\it Proof:} Routine using Lemma \ref{lem:scalar} and the meaning of inversion. \hfill $\Box$

\medskip
\noindent We now turn our attention to the inverted $\Phi^*$-split basis. 

\begin{lemma}
\label{lem:dualinvlbubbasis}
With reference to Notation \ref{not:aastar}, let $\lbrace v_i \rbrace_{i=0}^d$ denote any vectors in $V$. Then the following are equivalent.

\begin{enumerate}
 \item[\rm (i)] 
$\lbrace v_i \rbrace_{i=0}^d$ is an inverted $\Phi^*$-split basis for $V$.
 \item[\rm (ii)] 
There exists $0 \neq \eta^*_0 \in E^*_0V$ such that  
\begin{eqnarray*}
v_i = (A-\theta_{d-i+1} I)\cdots (A-\theta_{d-1} I)(A-\theta_{d}I)\eta^*_0  \qquad \qquad (0 \leq i \leq d).
\end{eqnarray*}
 \item[\rm (iii)] 
There exists $0 \neq \eta_0 \in E_0V$ such that  
\begin{eqnarray*}
\displaystyle{
v_i = \frac{(A^*-\theta^*_{i+1} I)\cdots (A^*-\theta^*_{d-1} I)(A^*-\theta^*_{d}I)\eta_0}{\phi_{i+1} \cdots \phi_{d-1} \phi_d}  \qquad \ \ (0 \leq i \leq d)},
\end{eqnarray*}
where $\lbrace \phi_i \rbrace_{i=1}^d$ is the split sequence of $\Phi$.
\end{enumerate}
\end{lemma}

\noindent {\it Proof:} Routine using Lemma \ref{lem:duallbubbasis} and the meaning of inversion. \hfill $\Box$

\begin{proposition}
\label{prop:dualinvlbublooklike}
With reference to Notation \ref{not:aastar}, let $B$ (resp. $B^*$) denote the matrix in $\hbox{Mat}_{d+1}(\K)$ which represents $A$ (resp. $A^*)$ with respect to an inverted $\Phi^*$-split basis for $V$. Then 
\begin{eqnarray*}
B = 
\left(
\begin{array}{c c c c c c}
\theta_d & & & & & {\bf 0} \\
1 & \theta_{d-1} &  & & & \\
& 1 & \theta_{d-2} &  & & \\
& & \cdot & \cdot &  &  \\
& & & \cdot & \cdot &  \\
{\bf 0}& & & & 1 & \theta_0
\end{array}
\right),
\qquad  \quad 
B^* = 
\left(
\begin{array}{c c c c c c}
\theta^*_0 & \phi_1 & & & & {\bf 0} \\
& \theta^*_1 & \phi_2 & & & \\
& & \theta^*_2 & \cdot & & \\
& &  & \cdot & \cdot &  \\
& & &  & \cdot & \phi_d \\
{\bf 0}& & & &  & \theta^*_d
\end{array}
\right),
\end{eqnarray*}
\noindent where $\lbrace \phi_i \rbrace_{i=1}^d$ is the split sequence of $\Phi$. In particular, $B$ is lower bidiagonal and $B^*$ is upper bidiagonal. 
\end{proposition}

\noindent {\it Proof:} Routine using Proposition \ref{prop:duallbublooklike} and the meaning of inversion. \hfill $\Box$

\begin{lemma}
With reference to Notation \ref{not:aastar}, let $\lbrace v_i \rbrace_{i=0}^d$ denote a $\Phi$-split basis for $V$, and let $\lbrace w_i \rbrace_{i=0}^d$ denote any vectors in $V$. Then the following are equivalent. 
\begin{enumerate}
\item[\rm (i)]  
$\lbrace w_i \rbrace_{i=0}^d$ is an inverted $\Phi^*$-split basis for $V$. 
\item[\rm (ii)] 
There exists $0 \neq c \in \K$ such that $w_i = c\,\phi_1 \phi_{2} \cdots \phi_{i}\, v_{i}$ for $0 \leq i \leq d$,
where $\lbrace \phi_i \rbrace_{i=1}^d$ is the split sequence of $\Phi$. 
\end{enumerate}
\end{lemma}

\noindent {\it Proof:} Routine using Lemma \ref{cor:trans} and the meaning of inversion.  \hfill $\Box$

\section{Classification of TH systems}
In this section we classify the TH systems up to isomorphism. 
\medskip

\noindent Let $\Phi$ denote a TH system on $V$. In the previous sections we associated with $\Phi$ some sequences of scalars: the eigenvalue sequence, the dual eigenvalue sequence, and the split sequence. We now show that those sequences determine $\Phi$ up to isomorphism.  

\begin{lemma}
\label{lem:paramsdetisoS99}
Let $\Phi$ and $\Phi'$ denote TH systems over $\fld$. Then the following are equivalent.
\begin{enumerate}
\item[\rm (i)]  $\Phi$ and $\Phi'$ are isomorphic.
\item[\rm (ii)]  $\Phi$ and $\Phi'$ share the same eigenvalue sequence, dual eigenvalue sequence, and split sequence. 
\end{enumerate}
\end{lemma}

\noindent {\it Proof:} (i) $\rightarrow$ (ii). Clear.

\noindent (ii) $\rightarrow$ (i). Without loss of generality, assume that $\Phi$ is the TH system from Definition \ref{def:HS}. Write $\Phi'=(A';\lbrace E_i' \rbrace_{i=0}^d; A^{*\prime};\lbrace E^{*\prime}_i\rbrace_{i=0}^d)$. Let $V'$ denote the vector space underlying $\Phi'$. Let $\lbrace v_i \rbrace_{i=0}^d$ (resp. $\lbrace v_i' \rbrace_{i=0}^d$) denote a $\Phi$-split (resp. $\Phi'$-split) basis for $V$ (resp. $V'$). Let $\gamma$ denote the vector space isomorphism from $V$ to $V'$ which sends $v_i$ to $v_i'$ for $0 \leq i \leq d$. Since $\Phi$ and $\Phi'$ share the same eigenvalue, dual eigenvalue, and split sequences, the matrix representing $A$ with respect to $\lbrace v_i \rbrace_{i=0}^d$ and the matrix representing $A'$ with respect to $\lbrace v_i' \rbrace_{i=0}^d$ are identical. Thus $\gamma A = A'\gamma$. Similarly $\gamma A^*= A^{*\prime}\gamma$. Moreover, using (\ref{eq:defEi}) we find $\gamma E_i = E_i' \gamma$ and $\gamma E^*_i = E^{*\prime}_i \gamma$ for  $0 \leq i \leq d$. Therefore $\gamma$ is an isomorphism of TH systems from $\Phi$ to $\Phi'$ and so $\Phi$ and $\Phi'$ are isomorphic. \hfill $\Box$

\bigskip

\noindent Lemma \ref{lem:paramsdetisoS99} motivates the following definition. 

\begin{definition}
\label{def:paofphi}
\rm
Let $\Phi$ denote a TH system on $V$. By the {\it parameter array of $\Phi$} we mean the sequence $(\lbrace \theta_i \rbrace_{i=0}^d, \lbrace \theta^*_i \rbrace_{i=0}^d, \lbrace \phi_i \rbrace_{i=1}^d)$, where 
$\lbrace \theta_i \rbrace_{i=0}^d$ (resp. $\lbrace \theta^*_i \rbrace_{i=0}^d$) is the eigenvalue (resp. dual eigenvalue) sequence of $\Phi$ and $\lbrace \phi_i \rbrace^d_{i=1}$ is the split sequence of $\Phi$. 
\end{definition}

\noindent In the following theorem, we classify the TH systems up to isomorphism. 

\begin{theorem}
\label{thm:classificationpa}
Let  
\begin{equation}
\label{eq:pa}
(\lbrace \theta_i \rbrace_{i=0}^d, \lbrace \theta^*_i \rbrace_{i=0}^d, \lbrace \phi_i \rbrace_{i=1}^d)
\end{equation}
denote scalars in $\K$. Then there exists a TH system $\Phi$ over $\K$ with parameter array (\ref{eq:pa}) if and only if {\rm (i)--(iii)} hold below. 

\begin{enumerate}
\item[\rm (i)] 
$ \theta_i\not=\theta_j \qquad $ \ if $\;\;i\not=j \qquad \qquad \qquad \ (0 \leq i,j\leq d)$.
\item[\rm (ii)]
$\theta^*_i\not=\theta^*_j \qquad $ \ \ if $\;\;i\not=j \qquad \qquad \qquad \ \ (0 \leq i,j\leq d)$.
\item[\rm (iii)] 
$\phi_i \not=0 \qquad \qquad \qquad \qquad \qquad \qquad \qquad (1 \leq i \leq d)$.
\end{enumerate}

\noindent Moreover if {\rm (i)--(iii)} hold above then $\Phi$ is unique up to isomorphism of TH systems.
\end{theorem}

\noindent {\it Proof:} To prove the theorem in one direction, let $\Phi=(A;\lbrace E_i\rbrace_{i=0}^d;A^*;\lbrace E^*_i\rbrace_{i=0}^d)$ denote a TH system over $\K$ with parameter array (\ref{eq:pa}). 
We show that conditions (i)--(iii) above hold. Conditions (i), (ii) hold by Definition \ref{def:evseq} and condition (iii) holds by the observation prior to Proposition \ref{prop:lbublooklike}. We are done with the proof in one direction.

\medskip

\noindent For the other direction, suppose conditions (i)--(iii) above hold. Consider the following matrices $B, B^* \in \hbox{Mat}_{d+1}(\K)$.
\beast
B = 
\left(
\begin{array}{c c c c c c}
\theta_d & & & & & {\bf 0} \\
\phi_1 & \theta_{d-1} &  & & & \\
& \phi_2 & \theta_{d-2} &  & & \\
& & \cdot & \cdot &  &  \\
& & & \cdot & \cdot &  \\
{\bf 0}& & & & \phi_d & \theta_0
\end{array}
\right),
&&\quad 
B^* = 
\left(
\begin{array}{c c c c c c}
\theta^*_0 & 1 & & & & {\bf 0} \\
& \theta^*_1 & 1 & & & \\
&  & \theta^*_2 & \cdot & & \\
& &  & \cdot & \cdot &  \\
& & &  & \cdot & 1 \\
{\bf 0}& & & &  & \theta^*_d
\end{array}
\right).
\eeast
Let $\lbrace u_i \rbrace_{i=0}^d$
denote a basis for $V$. Let $A$ (resp. $A^*$) denote the element of ${\rm End}(V)$ which is represented by $B$ (resp. $B^*$) with respect to $\lbrace u_i \rbrace_{i=0}^d$. We observe that $A$ (resp. $A^*$) is multiplicity-free, with eigenvalues $\lbrace \theta_i \rbrace_{i=0}^d$
(resp. $\lbrace \theta^*_i \rbrace_{i=0}^d$). For $0 \leq i \leq d$, let $E_i$ (resp. $E^*_i$) denote
the primitive idempotent of $A$ (resp. $A^*$) corresponding to $\theta_i$ (resp. $\theta^*_i$).
We show that $\Phi=(A;\lbrace E_i\rbrace_{i=0}^d;A^*;\lbrace E^*_i\rbrace_{i=0}^d)$ is a TH system on $V$.
To do this, we show that $\Phi$ satisfies conditions (i)--(v) in Definition \ref{def:HS}.
Conditions (i)--(iii) are clearly satisfied, so it remains to prove conditions (iv), (v). To prove condition (iv), we make a claim. For $0 \leq i \leq d$, let
$U_i = \hbox{Span}(u_i)$. For notational convenience, set $U_{-1}=0$ and $U_{d+1}=0$. Then 
\begin{eqnarray}
U_{i}+U_{i+1}+\cdots +U_{d} &=& E_0V + E_1V + \cdots + E_{d-i}V.
\label{eq:vsumitod}
\end{eqnarray}
\noindent To prove the claim, abbreviate $K=\sum^d_{h=i} U_{h}$ and $L=\sum^{d-i}_{h=0}E_hV$.
We show that $K=L$. To obtain $L\subseteq K$, set $X=\prod_{h=0}^{i-1} (A-\theta_{d-h}I)$, and observe that 
$L=XV$ by elementary linear algebra. Observe that $(A-\theta_{d-h}I)U_h = U_{h+1}$ for $0 \leq h \leq d$
and hence $XU_h\subseteq K$ for $0 \leq h \leq d$. Since $V = \sum_{h=0}^d U_h$ we have $XV\subseteq K$.
We now have $L\subseteq K$. To obtain $K\subseteq L$, set $Y=\prod_{h=i}^{d} (A-\theta_{d-h}I)$, and observe that
\begin{eqnarray}
L&=&\lbrace v \in V \;|\;Yv = 0\rbrace .
\label{eq:Keraction}
\end{eqnarray}
Since $(A-\theta_{d-h}I)U_h = U_{h+1}$ for $0 \leq h \leq d$, we have $YU_j=0$ for $ i \leq j \leq d$, so $YK=0$. Combining this with (\ref{eq:Keraction}), we find $K\subseteq L$. We now have $K=L$ and the claim is proved. We now prove condition (iv). 
Fix integers $i,j$ $(0 \leq i, j \leq d)$ such that $i-j \geq 1$. First assume that $i-j > 1$. Observe that $(A^*-\theta^*_hI)U_h = U_{h-1}$ for $0 \leq h \leq d$. Using this fact and  (\ref{eq:vsumitod}), we find $A^*E_{j}V \subseteq A^* \sum_{h=0}^{j} E_{h}V = A^* \sum_{h=d-j}^{d} U_{h} \subseteq \sum_{h=d-j-1}^{d} U_{h} = \sum_{h=0}^{j+1} E_{h}V$. Therefore $E_iA^*E_jV \subseteq E_i \sum_{h=0}^{j+1} E_{h}V = 0$ so $E_iA^*E_j = 0$.
Next assume that $i-j = 1$. We show that $E_{i}A^*E_j \neq 0$. By way of contradiction, assume that $E_{i}A^*E_j = 0$. Using these comments and (\ref{eq:vsumitod}), we find $A^* U_{d-j} \subseteq A^* \sum_{h=d-j}^{d} U_{h} = A^* \sum_{h=0}^{j} E_{h}V = (\sum^d_{k=0} E_k)A^* (\sum_{h=0}^{j} E_{h})V \subseteq \sum_{h=0}^{j} E_{h}V = \sum_{h=d-j}^{d} U_{h}$. This contradicts our earlier remark that $(A^*-\theta^*_{d-j}I)U_{d-j} = U_{d-j-1}$. Therefore $E_{i}A^*E_j \neq 0$ and Definition \ref{def:HS}(iv) holds. The proof for Definition \ref{def:HS}(v) is similar and omitted. We have now shown that $\Phi$ is a TH system on $V$. By Definition \ref{def:evseq}, $\Phi$ has eigenvalue sequence $\lbrace \theta_i \rbrace_{i=0}^d$ and dual eigenvalue sequence $\lbrace \theta^*_i \rbrace_{i=0}^d$. By Proposition \ref{thm:lbubbasisvsrep}, $\lbrace u_i \rbrace_{i=0}^d$ is a $\Phi$-split basis for $V$. Therefore $\Phi$ has split sequence $\lbrace \phi_i \rbrace_{i=1}^d$ by Proposition \ref{prop:lbublooklike}. The TH system $\Phi$ is unique up to isomorphism by Lemma \ref{lem:paramsdetisoS99}. \hfill $\Box$

\bigskip

\noindent We finish this section with a comment. 

\begin{lemma}
Let $\Phi$ denote a TH system on $V$ with parameter array $(\lbrace \theta_i \rbrace_{i=0}^d, \lbrace \theta^*_i \rbrace_{i=0}^d, \lbrace \phi_i \rbrace_{i=1}^d)$. Then the dual TH system $\Phi^*$ has parameter array $(\lbrace \theta^*_i \rbrace_{i=0}^d, \lbrace \theta_i \rbrace_{i=0}^d, \lbrace \phi_{d-i+1} \rbrace_{i=1}^d)$.
\end{lemma}

\noindent {\it Proof:} Immediate from Definitions \ref{def:THdual}, \ref{def:evseq} and Lemma \ref{lem:splitseqdual}. \hfill $\Box$

\section{The scalar $\nu$}
In this section we introduce a scalar $\nu$ that will help us describe TH systems. We start by updating Notation \ref{not:aastar}. 

\begin{notation}
\label{ass:updated}
\rm
Let $\Phi=(A;\lbrace E_i\rbrace_{i=0}^d;A^*;\lbrace E^*_i\rbrace_{i=0}^d)$ denote a TH system on $V$. 
Let \\ $(\lbrace \theta_i \rbrace_{i=0}^d, \lbrace \theta^*_i \rbrace_{i=0}^d, \lbrace \phi_i \rbrace_{i=1}^d)$ denote the parameter array of $\Phi$.
\end{notation}

\begin{lemma}
\label{cor:eta0gen}
With reference to Notation \ref{ass:updated}, let ${\mathcal D}$ (resp. ${\mathcal D}^*$) denote the $\fld$-subalgebra of ${\rm End}(V)$ generated by $A$ (resp. $A^*$). Fix $0 \neq \eta_0 \in E_0V$ and $0 \neq \eta^*_0 \in E^*_0V$. Then each of the maps
\begin{eqnarray*}
 {\mathcal D} \to V   \qquad \qquad \qquad \qquad \qquad  \; \; \, {\mathcal D}^* \to V \; \; \, \\
\; \, \, X \mapsto X\eta^*_0   \qquad \qquad \qquad \qquad \qquad  X \mapsto X\eta_0
\end{eqnarray*}
is an isomorphism of $\fld$-vector spaces. 
\end{lemma}

\noindent {\it Proof:} For the map on the left, use the fact that (\ref{eq:duallbubbasis}) is a basis for $V$.
For the map on the right, use the fact that (\ref{eq:spbasis}) is a basis for $V$. \hfill $\Box$

\begin{lemma}
\label{lem:eoeostar}
With reference to Notation \ref{ass:updated}, each of the maps
\begin{eqnarray*}
 E^*_0V &\to& E_0V    \qquad \qquad \qquad    \; \; \, E_0V \ \; \to \ \; E^*_0V  \; \; \, \\
v &\mapsto& E_0v    \qquad \qquad \qquad \qquad \; \;  v \ \; \mapsto \ \; E^*_0v
\end{eqnarray*}
is an isomorphism of $\fld$-vector spaces. 
\end{lemma}

\noindent {\it Proof:} Use Lemma \ref{cor:eta0gen} and the fact that each of $E_0V, E^*_0V$ has dimension $1$. \hfill $\Box$

\medskip

\noindent Consider the maps in Lemma \ref{lem:eoeostar}. If we compose the map on the left with the map on the right, then the resulting map acts on $E_0V$ as a nonzero scalar multiple of the identity. Denote the scalar by $\alpha$. If we now compose the map on the right with the map on the left, then the resulting map acts on $E^*_0V$ as $\alpha$ times the identity.  We define $\nu$ to be the reciprocal of $\alpha$. Observe that $\nu$ is nonzero. 

\begin{lemma}
\label{lem:nutripleproduct}
With reference to Notation \ref{ass:updated}, the following {\rm (i), (ii)} hold.
\begin{enumerate}
 \item[\rm (i)] 
$\nu E_0E^*_0E_0 = E_0.$
 \item[\rm (ii)]
$\nu E^*_0E_0E^*_0 = E^*_0.$
\end{enumerate}
\end{lemma}

\noindent {\it Proof:} Clear from the definition of $\nu$. \hfill $\Box$ 

\medskip

\noindent We mention one significance of $\nu$. 

\begin{lemma}
With reference to Notation \ref{ass:updated}, ${\rm tr}(E_0E^*_0) = \nu^{-1}$.
\end{lemma}

\noindent {\it Proof:} Consider the equation in Lemma \ref{lem:nutripleproduct}(i). Take the trace of each side and then simplify using the fact that ${\rm tr}(E_0) = 1$ and ${\rm tr}(E_0E^*_0E_0) = {\rm tr}(E_0E_0E^*_0) = {\rm tr}(E_0E^*_0)$. The result follows.  \hfill $\Box$

\medskip

\noindent We now express $\nu$ in terms of the parameter array of $\Phi$. 

\begin{lemma}
\label{lem:nupa}
With reference to Notation \ref{ass:updated}, 
\begin{eqnarray*}
\nu = \frac{(\theta_0-\theta_1)(\theta_0-\theta_2)\cdots(\theta_0-\theta_d)(\theta^*_0-\theta^*_1)(\theta^*_0-\theta^*_2)\cdots
(\theta^*_0-\theta^*_d)}
{\phi_1 \phi_2 \cdots \phi_d}.
\end{eqnarray*}
\end{lemma}

\noindent {\it Proof:} Fix $0 \neq \eta_0 \in E_0V$ and let $\lbrace v_i \rbrace_{i=0}^d$ denote the corresponding $\Phi$-split basis for $V$ from Definition \ref{def:lbubbasis}. Observe that $v_d = \eta_0$. 
By Lemma \ref{lem:nutripleproduct}, we have $\nu E_0E^*_0E_0 = E_0$. Applying each side of this equation to $v_d$ and using $v_d \in E_0V$, we find $\nu E_0E^*_0v_d = v_d$. By (\ref{eq:defEi}) we have $E_0 = \psi^{-1}\prod_{i=1}^d(A-\theta_iI)$, where $\psi = \prod_{i=1}^d(\theta_0-\theta_i)$. Similarly we have $E^*_0 = \psi^{*-1}\prod_{i=1}^d(A^*-\theta^*_iI)$, where $\psi^* = \prod_{i=1}^d (\theta^*_0-\theta^*_i)$. By Proposition \ref{prop:lbublooklike}, we find $\prod_{i=1}^d(A-\theta_iI)v_0 =  \phi \,v_d$, where $\phi=\prod_{i=1}^d \phi_i$. Similarly $\prod_{i=1}^d(A^*-\theta^*_iI)v_d = v_0$. Combining these facts, we find 
$\nu \phi v_d = \psi \psi^* v_d$. Now $\nu \phi = \psi \psi^*$ since $v_d \neq 0$. The result follows. \hfill $\Box$

\section{The $\Phi$-standard basis}
Let $\Phi$ denote a TH system on $V$. In this section we investigate a certain basis for $V$ called the $\Phi$-standard basis.

\begin{lemma}
With reference to Notation \ref{ass:updated}, let $0 \neq \eta_0 \in E_0V$. 
Then the sequence
\begin{eqnarray}
E^*_0\eta_0, E^*_1\eta_0, \ldots, E^*_d\eta_0
\label{eq:dst}
\end{eqnarray}
is a basis for $V$.
\end{lemma}

\noindent {\it Proof:} This follows from Lemma \ref{cor:eta0gen} and the fact that $\{E_i^*\}_{i=0}^d$ is a basis for ${\mathcal D}^*$.  \hfill $\Box$

\begin{definition} 
\label{def:hdbasis}
\rm
With reference to Notation \ref{ass:updated}, a basis for $V$ is called {\it $\Phi$-standard} whenever it is 
of the form (\ref{eq:dst}), where $0 \neq \eta_0 \in E_0V$.
\end{definition}

\begin{lemma}
With reference to Notation \ref{ass:updated}, let $\lbrace v_i \rbrace_{i=0}^d$ denote a $\Phi$-standard basis for $V$, and let $\lbrace w_i \rbrace_{i=0}^d$ denote any vectors in $V$. Then the following are equivalent. 

\begin{enumerate}
\item[\rm (i)]  
$\lbrace w_i \rbrace_{i=0}^d$ is a $\Phi$-standard basis for $V$. 
\item[\rm (ii)] 
There exists $0 \neq c \in \K$ such that $w_i = c \, v_i$ for $0 \leq i \leq d$.  
\end{enumerate}
\end{lemma}

\noindent {\it Proof:} Routine. \hfill $\Box$

\medskip

\begin{definition}
\label{not:ccstar}
\rm
With reference to Notation \ref{ass:updated}, let $H$ (resp. $D^*$) denote the matrix in $\hbox{Mat}_{d+1}(\K)$ which represents $A$ (resp. $A^*$) with respect to a $\Phi$-standard basis for $V$. 
\end{definition}

\noindent Our next goal is to describe the matrices $H$ and $D^*$. We start with $D^*$.  

\begin{lemma}
\label{lem:dstar}
With reference to Notation \ref{ass:updated} and Definition \ref{not:ccstar}, the matrix $D^*$ is diagonal with entries $D^*_{ii} = \theta^*_i$ for $0 \leq i \leq d$.
\end{lemma}

\noindent {\it Proof:} Recall $A^*E^*_i = \theta^*_i E^*_i$ for $0 \leq i \leq d$. The result follows.  \hfill $\Box$

\medskip
We now turn our attention to the matrix $H$. We recall some linear algebraic terms. Suppose we are given two bases for $V$, written $\lbrace u_i \rbrace_{i=0}^d$ and $\lbrace v_i \rbrace_{i=0}^d$. By the {\it transition matrix} from $\lbrace u_i \rbrace_{i=0}^d$ to $\lbrace v_i \rbrace_{i=0}^d$, we mean the matrix $T$ in $\hbox{Mat}_{d+1}(\K)$ that satisfies $v_j = \sum_{i=0}^d T_{ij}u_i$ for $0 \leq j\leq d$.
We recall a few properties of transition matrices. Let $T$ denote the transition matrix from $\lbrace u_i \rbrace_{i=0}^d$ to 
$\lbrace v_i \rbrace_{i=0}^d$. Then $T^{-1}$ exists, and equals the transition matrix from $\lbrace v_i \rbrace_{i=0}^d$ to
$\lbrace u_i \rbrace_{i=0}^d$. Let $\lbrace w_i \rbrace_{i=0}^d$ denote a basis for $V$, and let $S$ denote the transition matrix from $\lbrace v_i \rbrace_{i=0}^d$ to $\lbrace w_i \rbrace_{i=0}^d$. Then $TS$ is the transition matrix from $\lbrace u_i \rbrace_{i=0}^d$  to $\lbrace w_i \rbrace_{i=0}^d$. Let $A \in {\rm End}(V)$, and let $M$ (resp. $N$) denote the matrix in $\hbox{Mat}_{d+1}(\K)$ which represents $A$ with respect to $\lbrace u_i \rbrace_{i=0}^d$ (resp. $\lbrace v_i \rbrace_{i=0}^d$). Then $M = T N T^{-1}$.

\begin{lemma}
\label{lem:ttrans}
With reference to Notation \ref{ass:updated} and Definition \ref{not:ccstar}, let $\lbrace u_i \rbrace_{i=0}^d$ denote a $\Phi$-split basis for $V$, and let $\lbrace v_i \rbrace_{i=0}^d$ denote a $\Phi$-standard basis for $V$. 
Let $B$ (resp. $B^*$) denote the matrix in $\hbox{Mat}_{d+1}(\K)$ which represents $A$ (resp. $A^*$) with respect to $\lbrace u_i \rbrace_{i=0}^d$. 
Let $T \in \hbox{Mat}_{d+1}(\K)$ denote the transition matrix from $\lbrace v_i \rbrace_{i=0}^d$ to $\lbrace u_i \rbrace_{i=0}^d$. 
Then the following {\rm (i), (ii)} hold.
\begin{enumerate}
 \item[\rm (i)] 
$H = T B T^{-1}$.
 \item[\rm (ii)]
$D^* = T B^* T^{-1}$.
\end{enumerate}
\end{lemma}

\noindent {\it Proof:} Follows from the comments prior to this lemma. \hfill $\Box$

\medskip
Recall that we are trying to find the matrix $H$. To do this we use the equations in Lemma \ref{lem:ttrans}. In these equations the matrices $B$ and $B^*$ were found in Proposition \ref{prop:lbublooklike} and the matrix $D^*$ was found in Lemma \ref{lem:dstar}. We use Lemma \ref{lem:ttrans}(ii) to find the matrix $T$ and then use Lemma \ref{lem:ttrans}(i) to find the matrix $H$.  

\begin{lemma}
\label{lem:hdtolbubbasis}
With reference to Notation \ref{ass:updated}, fix $0 \neq \eta_0 \in E_0V$. Let $\lbrace u_i \rbrace_{i=0}^d$ denote the $\Phi$-split basis for $V$ as in (\ref{eq:spbasis}), and let $\lbrace v_i \rbrace_{i=0}^d$ denote the $\Phi$-standard basis for $V$ as in (\ref{eq:dst}). 
Let $T \in \hbox{Mat}_{d+1}(\K)$ denote the transition matrix from $\lbrace v_i \rbrace_{i=0}^d$ to $\lbrace u_i \rbrace_{i=0}^d$. Then $T$ is upper triangular with entries 
\begin{equation}
\label{eq:tij}
T_{ij} = (\theta^*_i - \theta^*_{j+1}) \cdots (\theta^*_i - \theta^*_{d-1})(\theta^*_i - \theta^*_d) \qquad \ \ (0 \leq i \leq j \leq d). 
\end{equation}
Moreover $T^{-1}$ is upper triangular with entries 
\begin{eqnarray*}
T^{-1}_{ij} = \frac{1}{(\theta^*_j - \theta^*_i)(\theta^*_j - \theta^*_{i+1}) \cdots (\theta^*_j - \theta^*_{j-1})} \, \frac{1}{(\theta^*_j - \theta^*_{j+1}) \cdots (\theta^*_j - \theta^*_{d-1})(\theta^*_j - \theta^*_d)} 
\end{eqnarray*}
for $0 \leq i \leq j \leq d$.
\end{lemma}

\noindent {\it Proof:} Observe that $u_d = \eta_0 = \sum_{i=0}^d E_i^* \eta_0 = \sum_{i=0}^d v_i$ and so $T_{id} = 1$ for $0 \leq i \leq d$. Moreover by Lemma \ref{lem:ttrans}(ii) we have $D^* T = T B^*$, where $D^*$ is from Lemma \ref{lem:dstar} and $B^*$ is from Proposition \ref{prop:lbublooklike}. Therefore by a routine matrix multiplication, we find $\theta^*_i T_{ij} = T_{i,j-1} + \theta^*_j T_{ij}$ for $0 \leq i \leq d$ and $1 \leq j \leq d$. Rearranging we obtain $T_{i,j-1} = (\theta^*_i - \theta^*_j) T_{ij}$ and (\ref{eq:tij}) follows by a simple recursion.
To prove our assertion about $T^{-1}$, let $S \in \hbox{Mat}_{d+1}(\K)$ denote an upper triangular matrix with entries 
\begin{equation}
\label{eq:sij}
S_{ij} = \frac{1}{(\theta^*_j - \theta^*_i)(\theta^*_j - \theta^*_{i+1}) \cdots (\theta^*_j - \theta^*_{j-1})} \, \frac{1}{(\theta^*_j - \theta^*_{j+1}) \cdots  (\theta^*_j - \theta^*_{d-1})(\theta^*_j - \theta^*_d)} 
\end{equation}
for $0 \leq i \leq j \leq d$.
It suffices to show that $TS = I$. The matrices $T$ and $S$ are both upper triangular, so $TS$
is upper triangular. By (\ref{eq:tij}), (\ref{eq:sij}) we find that for $0 \leq i \leq d$,
\beast
S_{ii} &=& \frac{1}{(\theta^*_i - \theta^*_{i+1}) \cdots (\theta^*_i - \theta^*_{d})}
\\
 &=& T_{ii}^{-1}
\eeast
so $(TS)_{ii} = 1$.
We now show that $(TS)_{ij}=0$ for $0 \leq i < j \leq d$. Let $i,j$ be given.  
It suffices to show that $(\theta^*_i-\theta^*_j)(TS)_{ij} = 0$, since $\lbrace \theta^*_h \rbrace_{h=0}^d$ are mutually distinct.
Observe that 
\begin{eqnarray}
(\theta^*_i-\theta^*_j)(TS)_{ij}
&=& (\theta^*_i-\theta^*_j)\sum_{h=0}^d T_{ih}S_{hj}    \nonumber
\\
&=& (\theta^*_i-\theta^*_j)\sum_{h=i}^j T_{ih}S_{hj}     \nonumber
\\
&=& \sum_{h=i}^j T_{ih}S_{hj}(\theta^*_i-\theta^*_h+\theta^*_h-\theta^*_j)
\nonumber
\\
&=& \sum_{h=i+1}^{j} T_{ih}(\theta^*_i-\theta^*_h)S_{hj}- 
\sum_{h=i}^{j-1} T_{ih}S_{hj}(\theta^*_j-\theta^*_h)
\nonumber
\\
&=& \sum_{h=i+1}^{j} T_{i,h-1}S_{hj}  -
\sum_{h=i}^{j-1} T_{ih}S_{h+1,j}
\label{eq:twosum}
\\
&=& 0
\nonumber
\end{eqnarray}
since the two sums in (\ref{eq:twosum}) are one and the same. We have now shown that $(TS)_{ij} = 0$ for $0 \leq i<j\leq d$. Combining our above arguments, we find $TS=I$. The result follows. \hfill $\Box$

\medskip

\begin{proposition}
\label{lem:hdlooklike}
With reference to Notation \ref{ass:updated} and Definition \ref{not:ccstar}, the matrix $H$ is Hessenberg with entries
\begin{eqnarray*}
\label{eq:hessentries}
H_{ij} &=& \sum^j_{h=i} \frac{(\theta^*_i - \theta^*_{h+1}) \cdots (\theta^*_i - \theta^*_{d-1})(\theta^*_i - \theta^*_{d}) \theta_{d-h}}{(\theta^*_j - \theta^*_{h})(\theta^*_j - \theta^*_{h+1}) \cdots (\theta^*_j - \theta^*_{j-1})(\theta^*_j - \theta^*_{j+1}) \cdots  (\theta^*_j - \theta^*_{d-1})(\theta^*_j - \theta^*_{d})} \nonumber 
\\
&+& \sum^{j- j^{\vee}}_{h= i- i'} \frac{(\theta^*_i - \theta^*_{h+2}) \cdots (\theta^*_i - \theta^*_{d-1})(\theta^*_i - \theta^*_{d}) \phi_{h+1}}{(\theta^*_j - \theta^*_{h})(\theta^*_j - \theta^*_{h+1}) \cdots (\theta^*_j - \theta^*_{j-1})(\theta^*_j - \theta^*_{j+1}) \cdots  (\theta^*_j - \theta^*_{d-1})(\theta^*_j - \theta^*_{d})},
\end{eqnarray*}
\noindent for $i \leq j+1 \; (0 \leq i,j \leq d)$, where ${\displaystyle{
i' = \cases{0, &if $\; i = 0$\cr
1, &if $\; i \neq 0$\cr}
}}$ and ${\displaystyle{
j^{\vee} = \cases{0, &if $\; j \neq d$\cr
1, &if $\; j = d.$\cr}}}$
\end{proposition}

\noindent {\it Proof:} By Lemma \ref{lem:ttrans}(i), $H = T B T^{-1}$, where $B$ is the matrix on the left in (\ref{eq:matrepaastar}) and $T, T^{-1}$ are from Lemma \ref{lem:hdtolbubbasis}. By a routine matrix multiplication, we find that each entry below the subdiagonal of $H$ is zero and the remaining entries of $H$ are as claimed. To see that $H$ is Hessenberg, it remains to show that each entry on the subdiagonal is nonzero. Using the above data, we find
$$H_{j+1,j} = \frac{(\theta^*_{j+1} - \theta^*_{j+2}) \cdots (\theta^*_{j+1} - \theta^*_{d}) \phi_{j+1}}{(\theta^*_j - \theta^*_{j+1}) \cdots (\theta^*_j - \theta^*_{d})},$$ for $0 \leq j \leq d-1$. Since $\lbrace \theta^*_i \rbrace_{i=0}^d$ are mutually distinct and $\lbrace \phi_i \rbrace_{i=1}^d$ are nonzero, $H_{j+1,j}$ is nonzero. Therefore $H$ is Hessenberg and the result follows. \hfill $\Box$

\bigskip
\noindent We give three characterizations of the $\Phi$-standard basis. 

\begin{proposition}
\label{lem:eggechar}
With reference to Notation \ref{ass:updated}, let $\lbrace v_i \rbrace_{i=0}^d$ denote a sequence of vectors in $V$, not all zero. Then this sequence is a $\Phi$-standard basis for $V$ if and only if the following {\rm (i), (ii)} hold.
\begin{enumerate}
\item[\rm (i)] $v_i \in E^*_iV$ for $0 \leq i \leq d$.
\item[\rm (ii)] $\sum_{i=0}^d v_i\in E_0V$.
\end{enumerate}
\end{proposition}

\noindent {\it Proof:} To prove the proposition in one direction, assume that $\lbrace v_i \rbrace_{i=0}^d$ is a $\Phi$-standard basis for $V$. By Definition \ref{def:hdbasis}, there exists $0 \neq \eta_0 \in E_0V$ such
that $v_i = E^*_i\eta_0$ for $0 \leq i \leq d$. Apparently $v_i \in E^*_iV$ so (i) holds. Observe that $I=\sum_{i=0}^d E^*_i$.
Applying both sides to $\eta_0$ we find $\eta_0=\sum_{i=0}^d v_i $ and (ii) follows. We have now proved the proposition in one direction. To prove the other direction, assume that $\lbrace v_i \rbrace_{i=0}^d$ satisfy (i), (ii) above. We define $\eta_0=\sum_{i=0}^d v_i$ and observe that $\eta_0 \in E_0V$ by (ii).
Using (i) we find $E^*_iv_j=\delta_{ij}v_j$ for $0 \leq i,j\leq d$ and hence $v_i = E^*_i\eta_0$ for $0 \leq i \leq d$. Observe that $\eta_0 \neq 0$ since at least one of $\lbrace v_i \rbrace_{i=0}^d$ is nonzero.
Now $\lbrace v_i \rbrace_{i=0}^d$ is a $\Phi$-standard basis for $V$ by Definition \ref{def:hdbasis}.
\hfill $\Box $\\

\noindent We recall some notation. For $X \in \hbox{Mat}_{d+1}(\K)$ and $\alpha \in \K$, $X$ is said to have {\it constant row sum $\alpha$} whenever $\sum_{h=0}^d X_{ih} = \alpha$ for $0 \leq i \leq d$.

\begin{proposition}
\label{lem:rowsum}
With reference to Notation \ref{ass:updated}, let $\lbrace v_i \rbrace_{i=0}^d$ denote a basis
for $V$. Let $C$ (resp. $C^*$) denote the matrix in $\hbox{Mat}_{d+1}(\K)$ which represents
$A$ (resp. $A^*)$ with respect to this basis. Then $\lbrace v_i \rbrace_{i=0}^d$ is a $\Phi$-standard basis for $V$ if and only if the following {\rm (i), (ii)} hold.
\begin{enumerate}
\item[\rm (i)] $C$ has constant row sum $\theta_0$.
\item[\rm (ii)] $C^*=\hbox{diag}(\theta^*_0, \theta^*_1, \ldots, \theta^*_d)$.
\end{enumerate}
\end{proposition}

\noindent {\it Proof:} Observe that $A \sum_{j=0}^d v_j = \sum_{i=0}^d v_i(C_{i0}+C_{i1}+\cdots+ C_{id})$.
Recall $E_0V$ is the eigenspace of $A$ corresponding to eigenvalue $\theta_0$.
Apparently $C$ has constant row sum $\theta_0$ if and only if $\sum_{i=0}^d v_i \in E_0V$.
Recall that for $0 \leq i \leq d$, $E^*_iV$ is the eigenspace of $A^*$ corresponding to eigenvalue $\theta^*_i$. Apparently $C^*=\hbox{diag}(\theta^*_0, \theta^*_1, \ldots, \theta^*_d)$
if and only if $v_i \in E^*_iV$ for $0 \leq i \leq d$. The result follows by Proposition
\ref{lem:eggechar}. \hfill $\Box $

\begin{proposition}
With reference to Notation \ref{ass:updated}, let $\lbrace v_i \rbrace_{i=0}^d$ denote a basis for $V$, and let $C$ (resp. $C^*$) denote the matrix in $\hbox{Mat}_{d+1}(\K)$ which represents $A$ (resp. $A^*$) with respect to this basis. Then $\lbrace v_i \rbrace_{i=0}^d$ is a $\Phi$-standard basis for $V$ if and only if the following {\rm (i), (ii)} hold.

\begin{enumerate}
\item[\rm (i)]
$C$ is Hessenberg with constant row sum $\theta_0$. 
\item[\rm (ii)]
$C^*$ is diagonal and $C^*_{00} = \theta^*_0$.
\end{enumerate}
\end{proposition}

\noindent {\it Proof:}
Suppose that $\lbrace v_i \rbrace_{i=0}^d$ is a $\Phi$-standard basis for $V$. By Proposition \ref{lem:hdlooklike}, $C$ is Hessenberg and by Proposition \ref{lem:rowsum}(i), $C$ has constant row sum $\theta_0$. By Proposition \ref{lem:rowsum}(ii), $C^*$ is diagonal and  $C^*_{00} = \theta^*_0$. Therefore (i), (ii) above hold. We have proved the proposition in one direction. To prove the other direction, suppose that (i), (ii) above hold.
Since $C$ is Hessenberg and $C^*$ is diagonal, applying Lemma \ref{lem:trihess} we find $\lbrace C^*_{ii} \rbrace_{i=0}^d$ is a dual eigenvalue sequence of $A, A^*$. Since $C^*_{00} = \theta^*_0$, applying Lemma \ref{lem:eig} we find $C^*_{ii} = \theta^*_i$ for $0 \leq i \leq d$. Therefore  $C^*=\hbox{diag}(\theta^*_0, \theta^*_1, \ldots, \theta^*_d)$ and by Proposition \ref{lem:rowsum}, $\lbrace v_i \rbrace_{i=0}^d$ is a $\Phi$-standard basis for $V$. \hfill $\Box$

\section{The $\Phi^*$-standard basis}
Let $\Phi=(A;\lbrace E_i\rbrace_{i=0}^d;A^*;\lbrace E^*_i\rbrace_{i=0}^d)$ denote a TH system on $V$. 
In the previous section, we obtained several results related to the $\Phi$-standard basis for $V$. Analogous results hold for the $\Phi^*$-standard basis for $V$. In this section, we display some of those results for use later in the paper. We begin by observing that a $\Phi^*$-standard basis for $V$ has the form
\begin{eqnarray}
E_0 \eta^*_0, E_1 \eta^*_0, \ldots, E_d \eta^*_0, \qquad \qquad \qquad 
\label{eq:st}
\end{eqnarray}
where $0 \neq \eta^*_0 \in E^*_0V$.

\begin{definition}
\label{not:ccstardual}
\rm
With reference to Notation \ref{ass:updated}, let $D$ (resp. $H^*$) denote the matrix in $\hbox{Mat}_{d+1}(\K)$ which represents $A$ (resp. $A^*$) with respect to a $\Phi^*$-standard basis for $V$. 
\end{definition}

\begin{lemma}
\label{lem:d}
With reference to Notation \ref{ass:updated} and Definition \ref{not:ccstardual}, the matrix $D$ is diagonal with entries $D_{ii} = \theta_i$ for $0 \leq i \leq d$.
\end{lemma}

\noindent {\it Proof:} Apply Lemma \ref{lem:dstar} to $\Phi^*$. \hfill $\Box$

\begin{lemma}
\label{lem:ttransdual}
With reference to Notation \ref{ass:updated} and Definition \ref{not:ccstardual}, let $\lbrace u_i \rbrace_{i=0}^d$ denote a $\Phi^*$-split basis for $V$, and let $\lbrace v_i \rbrace_{i=0}^d$ denote a $\Phi^*$-standard basis for $V$. Let $B$ (resp. $B^*$) denote the matrix in $\hbox{Mat}_{d+1}(\K)$ which represents $A$ (resp. $A^*$) with respect to $\lbrace u_i \rbrace_{i=0}^d$.      
Let $T^* \in \hbox{Mat}_{d+1}(\K)$ denote the transition matrix from $\lbrace v_i \rbrace_{i=0}^d$ to $\lbrace u_i \rbrace_{i=0}^d$. 
Then the following {\rm (i), (ii)} hold.
\begin{enumerate}
 \item[\rm (i)] 
$D = T^* B T^{*-1}$.
 \item[\rm (ii)]
$H^* = T^* B^* T^{*-1}$.
\end{enumerate}
\end{lemma}

\noindent {\it Proof:} Apply Lemma \ref{lem:ttrans} to $\Phi^*$. \hfill $\Box$

\begin{lemma}
\label{lem:dhtolbubbasis}
With reference to Notation \ref{ass:updated}, fix $0 \neq \eta^*_0 \in E^*_0V$. Let $\lbrace u_i \rbrace_{i=0}^d$ denote the $\Phi^*$-split basis for $V$ as in (\ref{eq:duallbubbasis}), and let $\lbrace v_i \rbrace_{i=0}^d$ denote the $\Phi^*$-standard basis for $V$ as in (\ref{eq:st}). 
Let $T^* \in \hbox{Mat}_{d+1}(\K)$ denote the transition matrix from $\lbrace v_i \rbrace_{i=0}^d$ to $\lbrace u_i \rbrace_{i=0}^d$. Then $T^*$ is upper triangular with entries 
\begin{eqnarray*}
T^*_{ij} = (\theta_i - \theta_{j+1}) \cdots (\theta_i - \theta_{d-1})(\theta_i - \theta_d) \qquad (0 \leq i \leq j \leq d). 
\end{eqnarray*}
Moreover $T^{*-1}$ is upper triangular with entries 
\begin{eqnarray*}
T^{*-1}_{ij} = \frac{1}{(\theta_j - \theta_i)(\theta_j - \theta_{i+1}) \cdots (\theta_j - \theta_{j-1})} \, \frac{1}{(\theta_j - \theta_{j+1}) \cdots (\theta_j - \theta_{d-1})(\theta_j - \theta_d)}
\end{eqnarray*}
for $0 \leq i \leq j \leq d$.
\end{lemma}

\noindent {\it Proof:} Apply Lemma \ref{lem:hdtolbubbasis} to $\Phi^*$. \hfill $\Box$

\begin{proposition}
\label{lem:hdlooklikedual}
With reference to Notation \ref{ass:updated} and Definition \ref{not:ccstardual}, the matrix
$H^*$ is Hessenberg with entries
\begin{eqnarray*}
H^*_{ij} &=& \sum^j_{h=i} \frac{(\theta_i - \theta_{h+1}) \cdots (\theta_i - \theta_{d-1})(\theta_i - \theta_{d}) \theta^*_{d-h}}{(\theta_j - \theta_{h})(\theta_j - \theta_{h+1}) \cdots (\theta_j - \theta_{j-1})(\theta_j - \theta_{j+1}) \cdots (\theta_j - \theta_{d-1})(\theta_j - \theta_{d})} \nonumber 
\\
&+& \sum^{j- j^{\vee}}_{h= i - i'} \frac{(\theta_i - \theta_{h+2}) \cdots (\theta_i - \theta_{d-1})(\theta_i - \theta_{d}) \phi_{d-h}}{(\theta_j - \theta_{h})(\theta_j - \theta_{h+1}) \cdots (\theta_j - \theta_{j-1})(\theta_j - \theta_{j+1}) \cdots (\theta_j - \theta_{d-1})(\theta_j - \theta_{d})},
\end{eqnarray*}
\noindent for $i \leq j+1 \; (0 \leq i,j \leq d)$, where ${\displaystyle{
i' = \cases{0, &if $\; i = 0$\cr
1, &if $\; i \neq 0$\cr}
}}$ and ${\displaystyle{
j^{\vee} = \cases{0, &if $\; j \neq d$\cr
1, &if $\; j = d.$\cr}}}$
\end{proposition}

\noindent {\it Proof:}  Apply Proposition \ref{lem:hdlooklike} to $\Phi^*$. \hfill $\Box$

\section{Transition matrices between the $\Phi$-standard basis and the $\Phi^*$-standard basis}
In this section we describe the transition matrices between the $\Phi$-standard basis for $V$ and the $\Phi^*$-standard basis for $V$. 

\begin{definition}
\label{def:zij}
\rm
With reference to Notation \ref{ass:updated}, let $Z$ denote the matrix in $\hbox{Mat}_{d+1}(\K)$ with entries  
\begin{eqnarray*}
\displaystyle{Z_{ij} = \delta_{i,d-j} \frac{(\theta_0 - \theta_1)(\theta_0 - \theta_2) \cdots (\theta_0 - \theta_d)}{\phi_{i+1} \phi_{i+2} \cdots \phi_d}} 
\end{eqnarray*}
for $0 \leq i, j \leq d$.
\end{definition}

\noindent The following lemma gives the significance of $Z$. 

\begin{lemma}
\label{lem:bdmatrix}
With reference to Notation \ref{ass:updated}, let $\lbrace v_i \rbrace_{i=0}^d$ denote a $\Phi$-split basis for $V$ as in (\ref{eq:spbasis}), and let $\lbrace w_i \rbrace_{i=0}^d$ denote a $\Phi^*$-split basis for $V$ as in (\ref{eq:duallbubbasis}) with $\eta_0, \eta^*_0$ from those lines chosen so that $\eta_0=E_0\eta^*_0$. Then the transition matrix from $\lbrace v_i \rbrace_{i=0}^d$ to $\lbrace w_i \rbrace_{i=0}^d$ is the matrix $Z$ from Definition \ref{def:zij}. 
\end{lemma}

\noindent {\it Proof:} The transition matrix from $\lbrace v_i \rbrace_{i=0}^d$ to $\lbrace w_i \rbrace_{i=0}^d$ is given by Lemma \ref{cor:trans} for an appropriate value of the scalar $c$ in that lemma. We now find that value. Setting $i=0$ in Lemma \ref{cor:trans}(ii), we find $w_0 = c\,v_d$. Using (\ref{eq:duallbubbasis}), we find $w_0 = (A - \theta_1I)(A - \theta_2I) \cdots (A - \theta_dI)\eta^*_0$ and this equals $(\theta_0 - \theta_1 )(\theta_0 - \theta_2 ) \cdots (\theta_0 - \theta_d ) E_0 \eta^*_0$. Using (\ref{eq:spbasis}), we find $v_d = \eta_0 = E_0 \eta^*_0$ by our choice of $\eta_0, \eta^*_0$. Therefore $w_0 = (\theta_0 - \theta_1 )(\theta_0 - \theta_2 ) \cdots (\theta_0 - \theta_d ) v_d$ and so $c = (\theta_0 - \theta_1)(\theta_0 - \theta_2) \cdots (\theta_0 - \theta_d)$. The result follows. \hfill $\Box$

\begin{definition}
\label{def:zijstar}
\rm
With reference to Notation \ref{ass:updated}, let $Z^*$ denote the matrix in $\hbox{Mat}_{d+1}(\K)$ with entries  
\begin{eqnarray*}
\displaystyle{Z^*_{ij} = \delta_{i,d-j} \frac{(\theta^*_0 - \theta^*_1)(\theta^*_0 - \theta^*_2) \cdots (\theta^*_0 - \theta^*_d)}{\phi_1 \phi_2 \cdots \phi_{d-i}}} 
\end{eqnarray*}
for $0 \leq i, j \leq d$.
\end{definition}

\noindent The following lemma gives the significance of $Z^*$. 

\begin{lemma}
\label{lem:bdmatrixstar}
With reference to Notation \ref{ass:updated}, let $\lbrace v_i \rbrace_{i=0}^d$ denote a $\Phi$-split basis for $V$ as in (\ref{eq:spbasis}), and let $\lbrace w_i \rbrace_{i=0}^d$ denote a $\Phi^*$-split basis for $V$ as in (\ref{eq:duallbubbasis}) with $\eta_0, \eta^*_0$ from those lines chosen so that $\eta^*_0=E^*_0\eta_0$. Then the transition matrix from $\lbrace w_i \rbrace_{i=0}^d$ to $\lbrace v_i \rbrace_{i=0}^d$ is the matrix $Z^*$ from Definition \ref{def:zijstar}. 
\end{lemma}

\noindent {\it Proof:} Apply Lemma \ref{lem:bdmatrix} to $\Phi^*$. \hfill $\Box$

\medskip

\begin{lemma}
\label{lem:nuzzstar}
With reference to Notation \ref{ass:updated} and Definitions \ref{def:zij}, \ref{def:zijstar}, 
\begin{eqnarray*}
ZZ^*=\nu I, \qquad \qquad Z^*Z = \nu I,
\end{eqnarray*}
where $\nu$ is the scalar above Lemma \ref{lem:nutripleproduct}.
\end{lemma}

\noindent {\it Proof:} Multiply $Z$ and $Z^*$ and use Lemma \ref{lem:nupa}. \hfill $\Box$

\begin{definition}
\label{def:ppstar}
\rm
With reference to Notation \ref{ass:updated}, let $P$ (resp. $P^*$) denote the transition matrix from (\ref{eq:st}) to (\ref{eq:dst})
(resp. (\ref{eq:dst}) to (\ref{eq:st})) with $\eta_0, \eta^*_0$ from those lines chosen so that $\eta^*_0=E^*_0\eta_0$ (resp. $\eta_0=E_0\eta^*_0$).
\end{definition}

\begin{lemma}
\label{lem:ptrans}
With reference to Notation \ref{ass:updated} and Definition \ref{def:ppstar}, 
\begin{eqnarray*}
P = T^* Z^* T^{-1}, \qquad \qquad P^* = T Z T^{*-1},
\end{eqnarray*}
where $T, T^*$ are from Lemmas \ref{lem:hdtolbubbasis}, \ref{lem:dhtolbubbasis} and $Z, Z^*$ are from Definitions \ref{def:zij}, \ref{def:zijstar}.  
\end{lemma}

\noindent {\it Proof:} By Lemma \ref{lem:hdtolbubbasis}, $T^{-1}$ is the transition matrix from (\ref{eq:spbasis}) to (\ref{eq:dst}). Moreover by Lemma \ref{lem:bdmatrixstar}, $Z^*$ is the transition matrix from (\ref{eq:duallbubbasis}) to (\ref{eq:spbasis}) and by Lemma \ref{lem:dhtolbubbasis}, $T^*$ is the transition matrix from (\ref{eq:st}) to (\ref{eq:duallbubbasis}). Therefore, by the comments prior to Lemma \ref{lem:ttrans}, we find $T^* Z^* T^{-1}$ is the transition matrix from (\ref{eq:st}) to (\ref{eq:dst}) and hence equals $P$. We have proved the equation on the left. To prove the equation on the right, apply the equation on the left to $\Phi^*$. \hfill $\Box$

\begin{theorem}
\label{thm:pentries}
With reference to Notation \ref{ass:updated} and Definition \ref{def:ppstar}, the following {\rm (i), (ii)} hold.
\begin{enumerate}
\item[\rm (i)]
\noindent For $0 \leq i,j \leq d$, $P_{ij}$ is equal to 
\begin{eqnarray*}
\frac{(\theta^*_0-\theta^*_1)(\theta^*_0-\theta^*_2)\cdots
(\theta^*_0-\theta^*_d)}
{(\theta^*_j-\theta^*_0)(\theta^*_j-\theta^*_1)\cdots (\theta^*_j-\theta^*_{j-1})
(\theta^*_j-\theta^*_{j+1}) \cdots (\theta^*_j - \theta^*_{d-1})(\theta^*_j-\theta^*_d)}
\end{eqnarray*}
times
\begin{eqnarray*}
\sum_{h=0}^d
\frac{(\theta_i-\theta_d)(\theta_i-\theta_{d-1})\cdots (\theta_i-\theta_{d-h+1})
(\theta^*_j-\theta^*_0)(\theta^*_j-\theta^*_1)
\cdots (\theta^*_j-\theta^*_{h-1})
}{\phi_1 \phi_2 \cdots \phi_h}.
\end{eqnarray*}
In particular $P_{i0}=1$ for $0 \leq i \leq d$. 
\item[\rm (ii)]
\noindent For $0 \leq i,j \leq d$, $P^*_{ij}$ is equal to 
\begin{eqnarray*}
\frac{(\theta_0-\theta_1)(\theta_0-\theta_2)\cdots
(\theta_0-\theta_d)}
{(\theta_j-\theta_0)(\theta_j-\theta_1)\cdots (\theta_j-\theta_{j-1})
(\theta_j-\theta_{j+1})\cdots (\theta_j - \theta_{d-1})(\theta_j-\theta_d)}
\end{eqnarray*}
times
\begin{eqnarray*}
\sum_{h=0}^d
\frac{(\theta_i^*-\theta_d^*)(\theta_i^*-\theta_{d-1}^*)\cdots (\theta_i^*-\theta_{d-h+1}^*)
(\theta_j-\theta_0)(\theta_j-\theta_1)
\cdots (\theta_j-\theta_{h-1})
}{\phi_d \phi_{d-1} \cdots \phi_{d-h+1}}.
\end{eqnarray*}
In particular $P^*_{i0}=1$ for $0 \leq i \leq d$.
\end{enumerate}
\end{theorem}

\noindent {\it Proof:} Routine by Lemma \ref{lem:ptrans} and matrix multiplication. \hfill $\Box$

\medskip

\begin{proposition}
With reference to Notation \ref{ass:updated} and Definition \ref{def:ppstar}, 
\begin{eqnarray*}
PP^*=\nu I, \qquad \qquad P^*P = \nu I,
\end{eqnarray*}
where $\nu$ is the scalar above Lemma \ref{lem:nutripleproduct}.
\end{proposition}

\noindent {\it Proof:} Evaluate the left hand sides using Lemma \ref{lem:ptrans} and then simplify using Lemma \ref{lem:nuzzstar}. \hfill $\Box$

\medskip

\noindent We finish this section with a comment. 

\begin{lemma}
With reference to Notation \ref{ass:updated} and Definition \ref{def:ppstar},   
\begin{eqnarray*}
H=P^{-1}DP, \qquad \qquad H^*=P^{*-1}D^*P^*,
\end{eqnarray*}
where $D^*$, $H$ (resp. $D$, $H^*$) are from Definition \ref{not:ccstar} (resp. Definition \ref{not:ccstardual}).
\end{lemma}

\noindent {\it Proof:} Recall from Definition \ref{not:ccstar} (resp. Definition \ref{not:ccstardual}) that $H$ (resp. $D$) is the matrix in $\hbox{Mat}_{d+1}(\K)$ which represents $A$ with respect to a $\Phi$-standard (resp. $\Phi^*$-standard) basis for $V$. Moreover recall from Definition \ref{def:ppstar} that $P$ is a transition matrix from a $\Phi^*$-standard basis to a $\Phi$-standard basis. Therefore the equation on the left holds by the comments prior to Lemma \ref{lem:ttrans}. To prove the equation on the right, apply the equation on the left to $\Phi^*$. \hfill $\Box$

\section{Acknowledgement}
This paper was written while the author was a graduate student at the University of Wisconsin-Madison. The author would like to thank his advisor Paul Terwilliger for his many valuable ideas and suggestions.

\bigskip

\noindent Ali Godjali \hfil\break
\noindent Department of Mathematics \hfil\break
\noindent University of Wisconsin \hfil\break
\noindent Van Vleck Hall \hfil\break
\noindent 480 Lincoln Drive \hfil\break
\noindent Madison, WI 53706-1388 USA \hfil\break
\noindent email: {\tt godjali@math.wisc.edu }\hfil\break

\end{document}